\documentclass[trsc,nonblindrev,oneside]{informs3} 
\OneAndAHalfSpacedXI 

\usepackage{graphicx}
\usepackage{algorithm}
\usepackage{algpseudocode}
\usepackage{booktabs}
\usepackage{optidef}
\usepackage{numprint}
\usepackage{ifdraft}
\usepackage[english]{babel}
\usepackage{pgf}
\usepackage{grffile}
\let\pgfimageWithoutPath\pgfimage
\renewcommand{\pgfimage}[2][]{\pgfimageWithoutPath[#1]{images/#2}}

\usepackage{caption}
\usepackage{subcaption}
\usepackage{glossaries}
\usepackage{csquotes}
\usepackage{tabularx}
\usepackage{tikz}

\usepackage{pgfplots}
\usepackage{tikzscale}
\usepackage[inline, shortlabels]{enumitem}
\setlist{itemjoin ={,\enspace},itemjoin* = { and\enspace}}

\pgfplotsset{compat=1.16}
\usepgfplotslibrary{groupplots}
\usetikzlibrary{graphs, graphs.standard, quotes,fit}

\usetikzlibrary{arrows.meta}
\usetikzlibrary{backgrounds}
\usepgfplotslibrary{patchplots}
\usepgfplotslibrary{fillbetween}

\usepackage{cleveref}

\usepackage{subfiles}

\definecolor{tum}{RGB}{0, 101, 189}

\definecolor{tum-darker}{RGB}{0, 51, 89}
\definecolor{tum-dark}{RGB}{0, 82, 147}
\definecolor{tum-medium}{RGB}{0, 115, 207}
\definecolor{tum-light}{RGB}{100, 160, 200}
\definecolor{tum-lighter}{RGB}{152, 198, 234}

\definecolor{tum-green}{RGB}{162, 173, 00}
\definecolor{tum-orange}{RGB}{227, 114, 34}
\definecolor{tum-ivory}{RGB}{218, 215, 203}

\newacronym{abk:ev}{EV}{electric vehicle}
\newacronym{abk:cpp}{CPP}{cutting plane problem}
\newacronym{abk:cgmp}{CGMP}{column generation master problem}
\newacronym{abk:soc}{SOC}{state-of-charge}
\newacronym{abk:osm}{OSM}{Open Street Maps}
\newacronym{abk:od}{OD}{Orgin-Destination}
\newacronym{abk:FRLM}{FRLM}{Flow Refueling Location Model}
\newacronym{abk:DFRLM}{DFRLM}{Deviation Flow Refueling Location Model}
\newacronym{abk:RSLP-R}{RSLP-R}{Refueling Station Location Problem
with Routing}
\newacronym{abk:SFCLM}{SFCLM}{Stochastic Flow-
Capturing Location Model}
\newacronym{abk:frlp}{FRLP}{Flow Refueling Location Problem}
\newacronym[
  longplural={Location Routing Problems with Intraroute Facilities}
]{abk:lrpif}{LRPIF}{Location Routing Problems with Intraroute Facilities}
\newacronym{abk:rcspp}{RCSPP}{Resource-Constrained Shortest Path Problem}
\newacronym{abk:REF}{REF}{resource extension function}
\newacronym{abk:VRP}{VRP}{Vehicle Routing Problem}

\newacronym{abk:ISA}{ISA}{Independent Scenario Approach}
\newacronym{abk:aSA}{$\alpha$-SA}{$\alpha$-Scenario Approach}
\newacronym{abk:aVA}{$\alpha$-VA}{$\alpha$-Vehicle Approach}
\newacronym{abk:FSA}{FSA}{Full Scenario Approach}

\usepackage{natbib}
 \bibpunct[, ]{(}{)}{,}{a}{}{,}%
 \def\bibfont{\small}%
\TheoremsNumberedThrough     

\EquationsNumberedThrough    

\MANUSCRIPTNO{TS-2022-0207} 

\begin{document}


\RUNAUTHOR{Godbersen, Kolisch, Schiffer}

\RUNTITLE{Robust Charging Station Network Planning}

\TITLE{Robust Charging Network Planning for Metropolitan Taxi Fleets}

\ARTICLEAUTHORS{%
\AUTHOR{Gregor Godbersen, Rainer Kolisch}
\AFF{School of Management, Department of Operations \& Technology \\ Technical University of Munich, Germany, \EMAIL{g.godbersen@tum.de}, \EMAIL{rainer.kolisch@tum.de}, \URL{https://www.tum.de}}

\AUTHOR{Maximilian Schiffer}
\AFF{School of Management \& Munich Data Science Institute \\
Technical University of Munich, Germany, \EMAIL{schiffer@tum.de}, \URL{https://www.tum.de}}

} 

\ABSTRACT{%
%
%

We study the robust charging station location problem for a large-scale commercial taxi fleet. 
Vehicles within the fleet coordinate on charging operations but not on customer acquisition. 
We decide on a set of charging stations to open to ensure operational feasibility. 
%
%
To take this decision, we propose a novel solution method situated between the Location Routing Problems with Intraroute Facilities and 
Flow Refueling Location Problems. Additionally, we introduce a problem variant that makes a station sizing decision. Using our exact approach, charging stations for a robust operation of city-wide taxi fleets can be planned.
%
%
We develop a deterministic core problem employing a cutting plane method for the strategic problem and a branch-and-price decomposition for the operational problem.
We embed this problem into a robust solution framework based on adversarial sampling, which allows for planner-selectable risk tolerance.
%
We solve instances derived from real-world data of the metropolitan area of Munich, containing 1,000 vehicles and 60 potential charging station locations.
Our investigation of the sensitivity of technological developments shows that increasing battery capacities show a more favorable impact on vehicle feasibility of up to 10 percentage points compared to increasing charging speeds. Allowing for depot charging dominates both of these options. Finally, we show that allowing just 1\% of operational infeasibility risk lowers infrastructure costs by 20\%.

}%


\HISTORY{First version: 17.06.2022, Revised: 14.4.2023}

\maketitle

\newcommand{\revnew}[1]{#1}

\section{Introduction}
\Glspl{abk:ev} are seen as a key element of future sustainable transportation systems as they allow for clean, low-emission transport, and may obtain even a zero-emission balance from a well-to-wheel perspective if powered by renewable energy sources. In urban environments, \glspl{abk:ev} offer additional advantages that may improve quality of life, e.g., due to reduced noise and fine dust emissions. However, despite increasing societal environmental awareness, governmental subsidies, and local low emission zone restrictions, the market uptake of \glspl{abk:ev} fell short of expectations \citep{agency_global_2019}. One root cause for this decelerated market uptake is the perceived range anxiety among private and commercial users. 

Focusing on private users, improving public charging infrastructure, as well as at-work and at-home charging opportunities begin to mitigate range anxiety concerns such that the market uptake of \glspl{abk:ev} encounters a slight acceleration. Focusing on commercial fleets, the diffusion of \glspl{abk:ev} remains still scarce. Here, companies are reluctant to rely on public charging infrastructure as they require permanent availability of charging infrastructure to preserve robust fleet operations.
Especially, commercial fleets with a high utilization, e.g., taxi and ride-hailing services require privately-owned charging infrastructure to preserve a sufficient service quality. This requirement increases an operator's investment costs, and vice versa, decreases its perception of economically viable deployment scenarios.

In this context, the optimal design of charging stations remains a core planning task to foster the adoption of \glspl{abk:ev}  for ride-hailing fleets. Optimally deploying the necessary charging stations requires a novel planning approach that goes beyond existing planning approaches by taking dedicated fleet operations and demand scenarios into account. In this context, ride-hailing fleet operators can influence the vehicles' recharging operations in between offering customer services via operator-to-driver communication or in case of (future) autonomous fleets via direct fleet control. With this additional decision layer, the resulting problem opens a new problem class in between flow-based charging infrastructure planning for private vehicles \citep[see, e.g.,][]{capar_arc_2013} and integrated location-routing based planning approaches for logistics fleets \citep[see, e.g.,][]{SchifferWalther2018}. While charging infrastructure design has been a vivid field of research for the above-mentioned applications (cf. Section~\ref{subsec:sota}), to the best of our knowledge no solution approach exists to determine optimal charging infrastructure for ride-hailing fleets, taking vehicle operations and recharging decisions into account. To close this gap, we develop a scenario-based robust optimization approach to determine optimal charging station infrastructure for  ride-hailing fleets.  In the following we first present the state of the art, before we detail our contributions and the organization of this paper.

\subsection{State of the Art}\label{subsec:sota}
Our planning problem relates to various research streams, ranging from \gls{abk:ev} operations and routing to charging infrastructure design. To keep this paper concise, we refer the interested reader to the survey of \cite{shen_optimization_2019} for a general overview on \gls{abk:ev} service operations, to the survey of \cite{kchaou-boujelben_charging_2021} for an overview of charging station location problems, and to the survey of \cite{SchifferSchneiderEtAl2018b} for an overview on electric vehicle routing problems. 
In the following, we focus on the fields of \glspl{abk:frlp} and \glspl{abk:lrpif}, which are most closely related to our planning problem. We also discuss approaches that are between the two fields in terms of how the coverage of demand is modelled.

Seminal papers on \glspl{abk:frlp} have been published by \cite{hodgson_flow-capturing_1990} and \cite{kuby_flow-refueling_2005}. 
They considered vehicles making round-trips between \gls{abk:od} pairs and aggregate these demands into flows which induce charging demand.
Any charging station placed along such an \gls{abk:od} flow can cover the charging demand. They seek the cost-minimal charging station configuration that covers all demand flows.
Several extensions of the \gls{abk:frlp} have been studied with a particular interest on electric vehicle charging infrastructure design:
\cite{KubyLim2007} specified the \gls{abk:frlp} for alternative fuel vehicles, while \cite{wang_locating_2009} introduced a set-covering formulation of the \gls{abk:frlp}, and \cite{upchurch_model_2009} studied a capacitated variant of the problem. \cite{wang_locating_2010} combined an \gls{abk:frlp} with a node-based component to account for intra-city demand. 

\cite{capar_arc_2013} and \cite{mirhassani_flexible_2012} proposed alternative problem formulations, which yield significantly improved computational times for exact solutions. \cite{mak_infrastructure_2013-1} focused on network design for battery swapping stations and considered uncertain demand. \cite{hosseini_refueling-station_2015} studied a two-stage planning approach that accounts for stationary and flexible (mobile) charging stations. \cite{wu_stochastic_2017-1} focused on an \gls{abk:frlp} with uncertain vehicle flows, solved by sample-average approximation. \cite{kchaou_boujelben_efficient_2019} modeled uncertain vehicle driving ranges while \cite{hosseini_robust_2021} considered both vehicle driving range and flow volume as uncertain parameters.

Other approaches focused explicitly on the representation of alternative routes and routing decisions. \cite{kim_deviation-flow_2012} proposed a problem variant that allows a deviation flow so that stations close to the shortest path can cover demand. \cite{li_heuristic_2014} introduced a multi-path model, which considers k-shortest path alternatives for every \gls{abk:od} pair. \cite{yildiz_branch_2016} extended potential detours on paths to cycles and used a path segment modeling approach to propose the first \gls{abk:frlp} with explicit routing decisions that can be solved with a branch-and-price algorithm. \cite{gopfert_branchcut_2019} and \cite{arslan_branch-and-cut_2019} proposed alternative branch-and-cut formulations to solve the \gls{abk:frlp} with routing decisions. 

Focusing on electrified logistics fleets, \glspl{abk:lrpif} have been studied to account for the interdependencies between vehicle routing and charging station siting decisions. \cite{YangSun2015} were the first to focus on the design of a battery swapping station network, while \cite{SchifferWalther2017} were the first to focus on charging stations. Further works in this field focused on efficient heuristic algorithms \citep{SchifferWalther2018}, dedicated real-world case studies \citep{schiffer_integrated_2021}, and robust optimization approaches \citep{SchifferWalther2018b}.

The demand coverage of some approaches can be placed between the full flexibility of location routing and the aggregated flow routing approaches. 
\cite{kang_strategic_2014} combined a set-covering model with a sequence of \gls{abk:od} trips obtained from household activity patterns.  \cite{tu_optimizing_2016} developed a spatial-temporal demand coverage approach using real-world taxi data, which they solve using a genetic algorithm. \cite{yildiz_urban_2019} extended the \gls{abk:frlp} with routing decisions to an urban setting with time-indexed station capacity restrictions and a stochastic demand model. \cite{brandstatter_location_2020} proposed a path-based charging station location problem for a car-sharing case where charging must be available close to where customers return vehicles.
\\

Concluding, existing \gls{abk:frlp} approaches lack the explicit consideration of a sequence of \gls{abk:od} trips allowing planning for single vehicles and extended fleet recharging decisions. 
Further, these approaches cannot be applied to our problem setting as they mostly focus on long-distance travel settings and consider inter-city \gls{abk:od} flows. To the best of our knowledge, \cite{yildiz_urban_2019} remains the only work that focuses on large-scale fleets in an urban setting. However, they obtain a tractable solution approach by using a simplified charging demand scheme. 
Existing \gls{abk:lrpif} approaches are not suitable for the problem setting of urban taxi fleets where less coordination is possible than in logistic fleets.
The possibility of creating entirely new tours by freely changing vehicle to customer assignments and replanning pickup times offered by the \gls{abk:lrpif} does not translate well to our problem setting, where coordination between vehicle operators is limited to charging operations instead of customer acquisition. 
Further, this added flexibility does lead to tractability issues when working with real-world data instances for city-wide taxi fleets.
Approaches between these two streams that consider paths and charging schedules for individual vehicles require high customization to the problem setting. Existing approaches have not sufficiently addressed the challenges posed by urban taxi fleets. 
As can be seen, no recent work in related research fields exists that could be straightforwardly extended to our problem setting.

\subsection{Contribution}
This paper aims to close the research gap outlined above by introducing the first methodological framework that allows designing optimal charging infrastructure for large ride-hailing fleets, taking charging station capacities, detours, and routing decisions for charging activities, partial recharges, as well as operational constraints and uncertain customer demand into account. Specifically, our contribution is two-fold:

First, we propose a novel exact solution approach for the robust location of electric charging stations for a private taxi or ride-hailing fleet. We model a setting where there is coordination between vehicles on charging decisions but not on customer acquisition and tour planning. \revnew{Specifically, we take a fixed customer schedule per vehicle and complement it with a charging schedule, including potential detours.} This places our approach in the middle of the existing \gls{abk:frlp} and \gls{abk:lrpif} problem formulations.   We consider several features found in a real-world problem setting, such as uncertain \revnew{customer schedules} based on real-world GPS data, a non-linear charging function, partial recharges, capacitated charging stations, and individual detours per vehicle within the solution approach. Nevertheless, we can solve a real-world-sized problem instance containing 1,000 vehicles and 60 potential charging station locations in less than 5 hours of computation time without aggregating customer demand. We develop two alternative algorithm variants within our scenario-based robust solution framework that incorporate a planner selectable risk tolerance. 

Second, we evaluate the solution approach using several instances sampled from real-world data. We assess the impact of changing technological parameters on the charging stations' layout. We consider battery and charging technology advancements and potential process improvements to allow depot charging. Our results show that depot charging and battery capacity have a stronger impact than charging speed. An increase in battery capacity shows an increase in vehicle feasibility up to 10 percentage points larger than when increasing charging speed. Further, we show that allowing for just 1\% of operational infeasibility risk lowers infrastructure costs by 20\%. Finally, we validate the robustness of our solutions under varying risk tolerances by employing an a-posteriori feasibility study. This feasibility study confirms that the vehicle feasibility predominantly stays within the selected risk tolerance even when considering an out-of-sample scenario set for the \revnew{vehicles' customer schedules}.
\\

The remainder of this paper is as follows. Section~\ref{sec:problem_description} states our problem definition, while Section~\ref{sec:solution_method} details our methodology. Section~\ref{sec:experimental_design} details our experimental design, while  Section~\ref{sec:results} shows computational results, and synthesizes managerial insights. Finally, Section~\ref{sec:conclusion} concludes this paper with a short summary.

\section{Problem Setting}\label{sec:problem_description}
In the following, we introduce our planning problem, which comprises a strategic and an operational decision from a fleet operator's perspective. At the strategic level, we decide on charging infrastructure investment decisions, i.e., we select a subset of charging stations from a set of possible candidate stations. At the operational level, we solve a vehicle scheduling problem, that integrates visits to charging stations into vehicle operations. Accordingly, we seek to find a cost-minimal charging infrastructure configuration that allows a taxi fleet to remain operational while serving all customer demand.

We formalize this problem setting as follows: let  ${\cal {G}} = ({\cal N},{\cal A})$ be a graph with node set ${\cal N}$ and arc set ${\cal A}$. The node set consists of two subsets, a set of candidate locations for charging stations ${\cal N}^s$ and a set of customer pickup and drop-off locations ${\cal N}^c$. Each arc $(i,j) \in {\cal A}$ denotes a route between two nodes, characterized by a distance $d_{ij}$, a travel time $\tau_{ij}$, and an energy consumption $e_{ij}$. Let~${\cal V}$ be the set of taxis which are operational during a discrete planning horizon ${\cal T}=\{0,\ldots,T\}$. Times~$t \in \mathcal{T}$ are equally spaced, and we refer to the span between two adjacent times $t$ and $t+1$ as a period.  We denote a vehicle's state of charge at time $t$ as $q_{v,t}$ with $q^\text{b}_v = q_{v,0}$ being its initial state of charge at $t = 0$. We use the charging function $\gamma(q,l)$ to denote the amount of energy charged, which depends on the vehicle's current state of charge $q$ and the charging duration $l$. 
Further, we associate with each vehicle $v\in {\cal V}$ an ordered set ${\cal C}_v$, which denotes a \revnew{schedule} of customer trips $c \in {\cal C}_v$ that must be served by vehicle $v$ in ${\cal T}$. We use a quadruple $(L^o_c, L^d_c, T^o_c, T^d_c)$ to characterize each customer trip $c$ based on its pickup location $L^o_c$, drop-off location $L^d_c$, start time at the pickup $T^o_c$, and arrival time at the drop-off $T^d_c$. With this notation, the problem's operational and strategic planning components hold as follows.

\paragraph{Operational Problem:}
At the operational level, we determine a feasible charging schedule for each vehicle $v$. Specifically, we complement the \revnew{vehicle's customer schedule} $\mathcal{C}_v$ by charging station visits between serving customers to keep the vehicle operational. Choosing the respective charging station visits entails decisions on the amount of energy recharged during each visit, i.e., each charging process's start and end time. We create such a charging schedule subject to the following constraints:
\begin{enumerate}
    \item Any charging station visit must be synchronized with providing service due to customer requests, i.e., a vehicle can only use the time in between customer trips for driving to and recharging at a charging station.
    \item A vehicle's state of charge must be within the interval $[q^{\min}, q^{\max}]$ during the whole planning period; moreover, it must be greater or equal to $q^\text{e}_v$ at the end of the planning horizon.
    \item At fleet level, capacity constraints of charging stations must not be violated, i.e., the maximum number of vehicles charging at a charging station in parallel is limited.
\end{enumerate}

\paragraph{Strategic Problem:}
At the strategic level, we decide on building charging infrastructure. Formally, we select a set of charging stations ${\cal S} \subseteq {\cal N}^s$ out of a set of candidate stations. We associate each station $s\in{\cal N}^s$ with a dedicated cost~\revnew{$o_s$} and a number of charging points $n_s$. We seek for a configuration ${\cal S}$ with the minimum cost such that all vehicles remain operational after solving the charge scheduling problem at the operational level. 
In practice,  \revnew{the vehicles' customer schedules, and thus operations vary} in between different periods. To account for such variation, we introduce a primitive uncertainty set ${\cal Z}$ that represents the uncertain \revnew{customer schedules}. A robust solution to the strategic problem must then ensure feasibility for all realizations of ${\cal Z}$.


\section{Methodology}\label{sec:solution_method}
In this section, we develop a scenario-based robust optimization approach to solve the planning problem defined in Section~\ref{sec:problem_description}. While we use scenarios to account for varying \revnew{customer schedules}, the core of our framework remains a deterministic algorithm that utilizes a cutting plane method in which we generate cuts by solving an operational problem via branch-and-price.
In the following, we first detail the deterministic problem, focusing on the  cutting plane method in Section~\ref{sec:solution_method:stategic_problem} and on its subproblem in Section~\ref{sec:operational_problem}.  We then show how to embed this algorithm into a scenario-based robust optimization framework in Section~\ref{sec:robust_master}. 
\revnew{Finally, we introduce a variant that makes a charging station sizing decision in Section~\ref{sec:variable_sizing}.}

We separate the problem into two stages: a strategic binary programming problem that decides on opening charging stations in the first stage and a second stage operational problem that ensures feasibility by discarding operationally infeasible solutions. This separation follows a rationale that is similar to the rationale of  the logic-based Benders formulations for 0-1 problems as discussed in \cite{hooker_logic-based_2003} or to the hybrid MILP and constraint programming approach introduced by~\cite{jain_algorithms_2001}. 

\subsection{\revnew{Cutting Plane Approach for the First Stage Problem}}\label{sec:solution_method:stategic_problem}

In the first stage, we use a binary program to obtain the minimal cost charging station configuration that is operationally feasible. We discard operationally infeasible solutions that enforce certain stations to be opened using an iterative cutting plane approach. 

We introduce binary variables $x_s\in\{1,0\}$ to indicate whether a station $s \in \mathcal{N}^s$ is opened ($x_s=1$) or closed ($x_s=0$),  and denote the set of feasibility cover cuts by \revnew{$\mathcal{J}$}, which is empty at the start of the cutting plane procedure. Then, our \gls{abk:cpp} reads:
\begin{mini!}
    {}{\sum_{s \in \mathcal{N}^s} \revnew{o_s} \cdot x_s\label{eq:benders_master:obj}}{\label{eq:benders_master}}{}
    \addConstraint{\sum_{s \in \mathcal{I}_j} x_s }{\geq 1}{\quad \revnew{\forall j \in \mathcal{J}}\label{eq:benders_master:cuts}}
    \addConstraint{ x_s                          }{\in \{0,1\},}{\quad \forall s \in \mathcal{N}^s.\label{eq:benders_master:domain}}
\end{mini!}
The objective function~\eqref{eq:benders_master:obj} minimizes the cost for opening charging stations. 
Constraints~\eqref{eq:benders_master:cuts} enforce for each cut $j \in {\cal J}$ that at least one station from the associated cover set ${\cal I}_j$ has to be opened, and Constraints~\eqref{eq:benders_master:domain} defines the domain of the variables $x_s$.

We then solve the problem as described in Algorithm~(\ref{alg:cmp}).
At every iteration, we solve the first-stage \gls{abk:cpp} to optimality (l.3), obtain the resulting solution $\overline{x}_s$, and subsequently check its operational feasibility using our second-stage operational subproblem (l.\ref{lst:line:test_oracle}).  If the subproblem is infeasible,  the set of currently closed charging stations $\mathcal{I} = \{ s \in \mathcal{N}^s\ |\ \overline{x}_s = 0 \}$ is a cover (l.\ref{lst:line:create_cover}). We then generate the cut \begin{equation}
     \sum_{s \in \mathcal{I}}{ x_s} \geq 1 
\end{equation} for cover $\mathcal{I}$, to cut off the current solution. We can reformulate the Cover Constraint \eqref{eq:benders_master:cuts}  to \begin{equation} 
    \sum_{s \in \mathcal{I}}{ (1-x_s)} \leq |\mathcal{I}| - 1, 
\end{equation} which constitutes a knapsack cover cut \citep[cf.][]{balas_facets_1975}.
These cuts lead to an operationally feasible and cost minimal charging station configuration since they fulfill the following two properties: i) a cut excludes at least the current solution, thus guaranteeing convergence in finite time, and ii) a cut does not exclude any feasible solutions, thus guaranteeing optimality as only infeasible solutions are excluded.

\begin{algorithm}
    \caption{Deterministic Cutting Plane Approach}\label{alg:cmp}
    \linespread{1}\selectfont
    
    \begin{algorithmic}[1]
    \State Initialize \gls{abk:cpp}\;\Comment{Section~\ref{sec:solution_method:stategic_problem}}
    \Repeat
    \State  $\overline{x}_s \gets $ Solve \gls{abk:cpp} and obtain optimal $x_s$ values\;
    \If{$\overline{x}_s$ is not operationally feasible}\label{lst:line:test_oracle}\Comment{Section~\ref{sec:operational_problem}}
        \State $\mathcal{I} \gets \{ s \in \mathcal{N}^s\ |\ \overline{x}_s = 0 \}$\label{lst:line:create_cover}\;
        \State $\mathcal{O} \gets $ \Call{Strengthen Cover}{$\mathcal{I}$}  \Comment{Algorithm~\ref{alg:lifting}} \;
        \ForAll{$\mathcal{P} \in \mathcal{O}$}
            \State Add constraint $\sum_{s \in \mathcal{P}} x_s \geq 1$ to \gls{abk:cpp}\;
        \EndFor
    \EndIf
     \Until{$\overline{x}_s$ is operationally feasible}
      \end{algorithmic}
    \end{algorithm}

\begin{proposition}
   The cuts $\sum_{s \in \mathcal{I}}{ x_s} \geq 1$ guarantee convergence for the \gls{abk:cpp}.
\end{proposition}
\proof{Proof\label{proof:cuts-convergence1}}
    Let  $\overline{x}_s$ be an operationally infeasible solution to the \gls{abk:cpp} with cover $\mathcal{I} = \{ s \in \mathcal{N}^s\ |\ \overline{x}_s = 0 \}$.
    Since all charging stations contained in the cover are closed  ($\sum_{s \in \mathcal{I}} \overline{x}_s = 0$), the derived cut $\sum_{s \in \mathcal{I}}{ x_s} \geq 1$ discards the solution. 
    If a feasible solution exists, the cut accordingly guarantees convergence in finite time because each cut discards a charging station configuration out of a finite number of total configurations.  $\hfill\blacksquare$
    \endproof{}
    \begin{proposition}
        The cuts $\sum_{s \in \mathcal{I}}{ x_s} \geq 1$ guarantee optimality for the \gls{abk:cpp}.
     \end{proposition}
    \proof{Proof\label{proof:cuts-convergence2}}
     Let $\overline{x}^1_s$ be a solution that was proven to be infeasible in a preceding iteration of Algorithm~\ref{alg:cmp}, leading to the cut $\sum_{s \in \mathcal{I}^1}{ x_s} \geq 1$ based on its set of closed charging stations $\mathcal{I}^1 = \{ s \in \mathcal{N}^s\ |\ \overline{x}^1_s = 0 \}$.  Now let $\overline{x}^2_s$ be any other solution found after $\overline{x}^1_s$ with $\mathcal{I}^2 = \{ s \in \mathcal{N}^s\ |\ \overline{x}^2_s = 0 \}$.
     We then distinguish two cases:
     i) If $\mathcal{I}^2 \not\supset \mathcal{I}^1$ holds, the cut does not discard solution $\overline{x}^2_s$ since at least one charging station is opened in solution $\overline{x}^2_s$ that was not opened in solution $\overline{x}^1_s$.
     i) otherwise, if $\mathcal{I}^2 \supset \mathcal{I}^1$ holds, $\overline{x}^2_s$ must also be infeasible as it contains more closed charging stations than $\overline{x}^1_s$. In this case, the cut  rightfully discards the solution $\overline{x}^2_s$. $\hfill\blacksquare$
\endproof{}

The effectiveness of this knapsack cut strongly depends on the size of the cover $\mathcal{I}$, and is initially very weak. We apply a cover strengthening heuristic (Algorithm~\ref{alg:lifting}) to generate a set $\mathcal{O}$ of smaller and thus tighter covers from $\mathcal{I}$.
Given a cover $\mathcal{I}$, we choose a random-sized sample of the locations in $\mathcal{I}$ to find a smaller cover $\mathcal{\tilde{I}} \subset\mathcal{I}$ (l.\ref{lst:line:random_subset}). 
We then solve the second stage operational problem with all stations $s \in \mathcal{\tilde{I}}$ closed and all stations $s \in N^s \setminus \mathcal{\tilde{I}}$ opened. If this problem is infeasible, $\mathcal{\tilde{I}}$ is a valid cover, and we store $\mathcal{\tilde{I}}$ in the  set $\mathcal{O}$ (l.\ref{lst:line:add_to_O}). We continue to find smaller covers from $\mathcal{\tilde{I}}$.
 If the operational problem is feasible for some candidate cover $\mathcal{\tilde{I}}$, the latter cannot be added to the set of covers, and we backtrack the search to the previous cover (l.\ref{lst:line:break_on_feasible}).

Once the heuristic terminates, we prune the  set $\mathcal{O}$ to retain only those covers that are non-dominated, i.e., we remove a cover if it is a superset of another cover. We then add all non-dominated covers in $\mathcal{O}$ to  \revnew{$\mathcal{J}$}, resulting in multiple new cuts in the \gls{abk:cpp}. Since every new cover in \revnew{$\mathcal{J}$} must dominate the original cover $\mathcal{I}$, we still cut off the current solution while strongly improving the convergence of the \gls{abk:cpp}. 
Since the strengthening of covers requires repeated evaluation of the operational problem and thus can come at a considerable cost, we balance the time between searching for tighter cuts and applying the cuts generated so far. To do so, we set a time limit for the \revnew{strengthening} heuristic proportional to the time spent in the last iteration of the \gls{abk:cpp} (l.\ref{lst:line:time_limit}).

\begin{algorithm}
\caption{Cover strengthening heuristic}\label{alg:lifting}
\linespread{1}\selectfont

\begin{algorithmic}[1]
 \Function{Strengthen Cover}{$\mathcal{I}$} \label{lst:line:cover_lift}
\State $\mathcal{O} \gets \{\mathcal{I}\}$\;
 \While{time limit not exceeded}\label{lst:line:time_limit}
   \State $\mathcal{\tilde{I}} \gets \mathcal{I}$
    \While{$\mathcal{\tilde{I}} \neq \emptyset$} 
        \State $\mathcal{\tilde{I}} \gets$ Choose random subset of $\mathcal{\tilde{I}}$\label{lst:line:random_subset}
        \If{Closing all stations in $\mathcal{\tilde{I}}$ is not operationally feasible}
            \State $\mathcal{O} \gets \mathcal{O} \cup \{\mathcal{\tilde{I}}\}$\label{lst:line:add_to_O}
        \Else
            \State break\label{lst:line:break_on_feasible}
        \EndIf
    \EndWhile
   \EndWhile
   \State Prune $\mathcal{O}$ to retain only non-dominated covers\;
   \State \Return $\mathcal{O}$
 \EndFunction
  \end{algorithmic}
\end{algorithm}

\subsection{\revnew{Branch-and-Price Approach for the Second Stage Problem}}\label{sec:operational_problem}
In the second-stage operational problem, we search for a feasible vehicle schedule for a charging station configuration $\overline{x}_s$ obtained from the first-stage  \gls{abk:cpp}. This problem shares similarities with the fixed route vehicle charging problem introduced by \cite{montoya_electric_2017}. For each vehicle, we determine a route in which it visits charging stations between customer service trips to keep the vehicle operational without violating charging station capacity restrictions. 
We note that we use the  operational problem to evaluate the feasibility of a charging station configuration such that we need to prove feasibility but not optimality.

We decompose the operational problem by using branch-and-price and derive a column generation master problem and  its pricing subproblem. In the following, we detail the  column generation master problem~(\ref{subsubsec:cgmp}), its pricing subproblem~(\ref{subsubsec:cgpp}), and our branching strategies~(\ref{subsubsec:cgbs}).

\subsubsection{Column Generation Master Problem}\label{subsubsec:cgmp}
Given a set~$\mathcal{P}_v$ that denotes all time and \gls{abk:soc} feasible routes for vehicle~$v\in \mathcal{V}$ and a dummy route~$0$, our \gls{abk:cgmp} selects a route $p\in\mathcal{P}_v\cup\{0\}$ for each vehicle $v$, such that the number of vehicles charging at each charging station does not exceed the station's capacity at any time $t\in\mathcal{T}$. For each vehicle $v$ the dummy route~$0$  can be selected if no feasible, conflict free route in $\mathcal{P}_v$ can be chosen. We use a binary parameter~$k_{vpst}$ to denote whether route $p\in\mathcal{P}_v$ of vehicle $v$ occupies a charging station~$s$  during time span $[t,t+1)$ ($k_{vpst}=1$) or not ($k_{vpst}=0$). To select a route for vehicle~$v$, we use binary variables~$y_{vp}$, which indicates if route $p\in P_v\cup\{0\}$ of vehicle $v$ is selected ($y_{vp} = 1$) or not ($y_{vp} =0$).  
Furthermore, $n_s$~denotes the \revnew{number of charging points of a charging station $s$, i.e, its capacity,} and  $\overline{x}_s \in \{1,0\}$ denotes if station~$s$ has been selected ($\overline{x}_s=1$) in the \gls{abk:cpp} \eqref{eq:benders_master} or not ($\overline{x}_s=0$). With this notation, the \gls{abk:cgmp} is as follows: 
\begin{mini!}
    {}{ \sum_{v \in \mathcal{V}} y_{v0}}{\label{eq:colcg_master}}{\label{eq:colcg_master_obj}}
    \addConstraint{\sum_{p \in \mathcal{P}_v\cup\{0\}} y_{vp}}{=1}{\quad \forall v \in \mathcal{V},\quad (\rho_v)\label{eq:colcg_master_convexity}}
    \addConstraint{ \sum_{v \in \mathcal{V}}\sum_{p \in \mathcal{P}_v}  k_{vpst} \cdot y_{vp}  }{\leq n_s \cdot \overline{x}_s ,}{\quad \forall s \in \mathcal{N}^s, \forall t \in T,  \quad (\lambda_{st})\label{eq:colcg_master_capacity}}
    \addConstraint{ y_{vp}}{\in \{0,1\},}{\quad  \forall v \in \mathcal{V},\forall p \in \mathcal{P}_v\cup\{0\}}\label{eq:colcg_master_domain}
\end{mini!}
Objective~\eqref{eq:colcg_master_obj} minimizes the number of vehicles for which the dummy route has to be chosen if no other conflict-free route can be found.
 Accordingly, if the objective function of an optimal operational problem solution is larger than zero, the associated solution of the master problem is not feasible.
Vehicle convexity constraints~\eqref{eq:colcg_master_convexity} enforce the selection of a single route per vehicle. Constraints~\eqref{eq:colcg_master_capacity} enforce capacity constraints at each charging station and point in time. Finally, Constraints~\eqref{eq:colcg_master_domain} state the variable domains. 

We use a column-generation approach to generate promising routes $p\in\mathcal{P}_v$. Accordingly, we relax the integrality constraints \eqref{eq:colcg_master_domain}, denote the dual value associated with
constraints~(\ref{eq:colcg_master_convexity}) and~(\ref{eq:colcg_master_capacity}) 
as $\rho_v$ and $\lambda_{st}$ respectively, and derive the reduced cost of $y_{vp}$ as $\revnew{\beta_{vp}} = -  \rho_v  - 
\sum_{s \in \mathcal{N}^s, t \in \mathcal{T}} k_{vpst} \cdot \lambda_{st}$.
At every iteration of the column generation, we use this pricing information to
determine all new columns with negative reduced costs and add them to the \gls{abk:cgmp}.

\subsubsection{Pricing Subproblem }\label{subsubsec:cgpp}
The~\gls{abk:cgmp} shows a distinct block structure for each vehicle, in which only Constraints~\eqref{eq:colcg_master_capacity} couple routes of different vehicles. Accordingly, we solve a separate pricing subproblem, which we model as a \gls{abk:rcspp}, for each vehicle.

For each vehicle $v$, we derive a time-expanded network  $\mathcal{G}_v = (\mathcal{V}_v,\mathcal{A}_v)$ from the physical network $\mathcal{G}$, limiting the customer vertices ${\cal N}^c$ to those of vehicle $v$ and using the set of discrete times $\mathcal{T}$ for the time expansion.  This network representation abstracts possible vehicle operations from the underlying road graph. Its vertex set $\mathcal{V}_v = \mathcal{H}^o_v \cup  \mathcal{H}^d_v \cup \mathcal{Q}_v$ consists of three subsets:
\begin{enumerate}
    \item $\mathcal{H}^o_v$, containing a start vertex for each customer trip $c$ at time $T^o_c$
    \item $\mathcal{H}^d_v$, containing an end vertex for each customer trip $c$ at time $T^d_c$
    \item  $\mathcal{Q}_v$, containing charging vertices between directly succeeding customer trips for each feasible charging start time $t$ at charging station $s \in \mathcal{S}$. We consider only those start times to be feasible, that allow for at least one period of charging between the two customers.  For two succeeding customers $c$ and $c+1$ the set of charging vertices $\mathcal{Q}_v$ results to   $\mathcal{Q}_v = \{(s,t)\ \forall s\in\mathcal{S}, t\in\mathcal{T},c\in\mathcal{C}_v\ |\ T^d_c + \tau_{L^d_c,s} \leq t \leq T^o_{c+1}- \tau_{s,L^o_{c+1}} - 1\}$.
\end{enumerate}

\noindent  The arc set $\mathcal{A}_v = \mathcal{R}_v \cup \mathcal{P}_v \cup \mathcal{W}_v \cup \mathcal{B}^1_v  \cup \mathcal{B}^2_v $ consists of five subsets:
\begin{enumerate}
    \item $\mathcal{R}_v$,  containing an arc from its associated start vertex in $\mathcal{H}^o_v$ to its associated end vertex in  $\mathcal{H}^d_v$  for each customer $c$.
    \item  $\mathcal{P}_v$, containing a transfer arc from the end vertex of customer $c$ to the start vertex of customer $c+1$ for each sequence of customers $c$ and $c+1$. When using this arc, a vehicle moves directly from the drop-of destination of customer $c$ to the pick-up destination of customer $c+1$ without charging.
    \item $\mathcal{B}^1_v$, containing a traveling arc from the end vertex of customer $c$ to each feasible charging vertex for each sequence of customers $c$ and $c+1$.
    \item $\mathcal{B}^2_v$, containing arcs from each feasible charging vertex to the start vertex of customer $c+1$ for each sequence of customers $c$ and $c+1$.

    \item $\mathcal{W}_v$, containing holdover arcs between each pair of charging vertices in time $t$ and the succeeding time $t+1$.
\end{enumerate}

\begin{example}[Time-Expanded Network]
    Figure~\ref{fig:time_space_network} shows an example segment of such a time-expanded network with three customers and three charging stations. From left to right, it consists of the end vertex of customer $c_0$ in the vertex set $\mathcal{H}^d_v$, six charging vertices in $\mathcal{Q}_v$ for the three charging stations and their feasible charging periods, followed by the start vertex of customer $c_1$ in $\mathcal{H}^o_v$, the end vertex of customer $c_1$ in $\mathcal{H}^d_v$, three charging  vertices in $\mathcal{Q}_v$ for charging stations $s_1$ and $s_3$,  and finally the start vertex of customer $c_2$ in  $\mathcal{H}^o_v$.

    \begin{figure}[!h]
        \FIGURE
        {\ifdraft{}{\includegraphics[width=1\linewidth]{images/time_space_network}}}
        {Example of a time-expanded network for two idle times, three customer trips, and three charging stations.\label{fig:time_space_network}}
        {}
    \end{figure}
\end{example}

As vertex attributes, we introduce a binary charge indicator $z_{i}$ and charging costs $\pi_{i}$. The charge indicator shows if the battery is charged in vertex~$i$ $(\pi_{i}=1)$ or not $(\pi_{i}=0)$, while the charging costs reflect the pricing information obtained from the $\gls{abk:cgmp}$. For all charging vertices $k\in\mathcal{Q}_v$, we set the charging indicator $z_{k}$ to one and set $\pi_{k}$ to the $\gls{abk:cgmp}['s]$ dual value $\lambda_{st}$ of the site $s$ and time $t$ associated with charging vertex~$k$. We set $\pi_{i}$ and $z_{i}$ of the remaining vertices~$i \in \mathcal{V}_v\setminus\mathcal{Q}_v$ to zero, as we do not perform any charging operations when visiting them.
As arc attributes, we introduce range costs $r_{ij}$ reflecting the \gls{abk:soc} reduction incurred by traversing an arc.
For arcs $(i,j)\in \mathcal{R}_v  \cup \mathcal{B}^1_v  \cup \mathcal{B}^2_v  \cup \mathcal{P}_v$, we set  $r_{ij}$  to the energy consumption $e_{kl}$  of their corresponding arcs $(k,l)\in\mathcal{A}$ in the original non-expanded graph $G$. Since the vehicle remains stationary during holdover arcs $(i,j)\in\mathcal{W}_v$, we set their range cost to zero.

We solve the \gls{abk:rcspp} by applying a dynamic programming based label setting algorithm~\citep[cf.][]{irnich_shortest_2005} and represent a partial path from the first segment's start node to a vertex~$i$ by a label $L_i = ( C_i, R_i)$, where $C_i$ denotes the cost  accumulated along the path to vertex~$i$  and $R_i$ denotes the \gls{abk:soc} at vertex~$i$.
We initialize the first label of a partial path with $C_0 = 0$ and $R_0=q^b_v$, and extend label $L_i$ along
an arc $(i,j) \in A_v$ using the following \glspl{abk:REF}
\begin{align}
    C_j &= C_i + \pi_{i} \label{eq:refC} \\
    R_j &=  - r_{ij} + \begin{cases}
          R_i , & \text{if}\ z_i=0 \\
          \min(q^{\max}, R_i + \gamma(R_i)), & \text{if}\ z_i=1 \\
    \end{cases}\label{eq:refR}
\end{align}
We add the charging cost~$\pi_{i}$ to the accumulated charging costs in Equation~\eqref{eq:refC}. For the update of the  \gls{abk:soc} in \eqref{eq:refR}, we first subtract the range cost~$r_{ij}$ and then distinguish two cases based on the charging indicator~$z_i$ of vertex~$i$. If charging takes place $(z_i=1)$, we determine the amount charged through the charging function ~$\gamma(R_i)$ and limit the total amount charged to the maximum battery capacity~$q^{\max}$; otherwise, we keep the previous \gls{abk:soc}~$R_i$.
Note that we cannot precalculate the amount charged as the charging function is dependent on the previous \gls{abk:soc}~$R_i$, which is only known after evaluating the \gls{abk:REF}. We add the amount charged in vertex~$i$ to the \gls{abk:soc} when extending an outgoing arc $(i,j)$, since the charged amount is only available at the end of the period associated with vertex $i$. By doing so, we ensure that the charged amount is only available in labels extended along outgoing arcs $(i,j)\in\mathcal{B}^2_v$, which reach the next customer trip.

 We ensure feasibility for all labels by discarding any extended label $L_j$ if the minimum SOC requirement $R_j \geq q^{\min}_v$ is violated. Any labels extended to the last segment's destination must additionally fulfill the end-of-shift SOC requirement $R_{\text{last}} \geq q^\text{e}_v$. In addition to these feasibility checks, we also apply the following dominance rule to reduce the size of the search space.
\begin{definition}[Label Domination]
Let $L_k = (C_k, R_k),\ k \in \{a,b\}$ be two labels associated with
different paths ending in the same vertex. Label $L_a$ dominates label $L_b$ if, $C_a \leq C_b \wedge R_a \geq R_b$. Here, we ensure by $R_a \geq R_b$ that any future extension on label~$L_b$ can also be reached by extending label~$L_a$, while $C_a \leq C_b$  ensures that its cost is lower or equal.
\end{definition}
If $L_a$ dominates $L_b$, we discard $L_b$ accordingly.

Once every label has been extended, labels associated with the final node of the final segment in the shift describe
all non-dominated operationally feasible vehicle schedules. Each of these labels describes a path $p$, whose sequence of visited vertices $U_{p}$ can be determined by backtracking.
Since the label attribute $C_i$ represents the sum of  $\lambda_{st}$ values for all visited charging vertices, we can then compute the reduced costs of a path $p$ of vehicle $v$ as $c_{vp} = -  \rho_v  - C_i$.
We add all paths with negative reduced costs ($c_{vp} < 0$) to the \gls{abk:cgmp}.

\subsubsection{Branching}\label{subsubsec:cgbs}
If the optimum solution of the relaxed column generation master problem (\ref{eq:colcg_master}a-\ref{eq:colcg_master}d) is fractional, there has to be at least one vehicle $v\in\mathcal{V}$ for which the set of selected paths $\mathcal{D}_v = \{ p \in \mathcal{P}_v\ |\ y_{vp} > 0 \}$ contains more than one path, and thus the paths' sequences of visited vertices $U_{p}$ must deviate in at least one vertex~$i\in \mathcal{V}_v$, i.e., $ U_p \neq U_q,\ p,q \in \mathcal{D}_v,\ p \neq q$. Our branching strategy separates the problem on these vertices.  We branch on the arc variables of the subproblem to obtain an integer solution \citep[cf.][]{desrosiers_chapter_1995}, which allows us to embed the resulting constraints in the time-expanded network's graph structure. Specifically, we enforce visiting specific vertices by removing arcs leading to other alternatives and forbid visiting vertices by removing incoming arcs. 

\begin{definition}[Station and time deviation]
Given two paths $\{p,q\}$, we check for a deviation by pairwise comparison between the sequences of visited vertices $U_{p}$ and $U_{q}$.
Let $V_a = (s_a,t_a)$ and $V_b = (s_b,t_b)$ be two vertices. 
For $s_a \neq s_b$, we define this as a station deviation, for $s_a = s_b \land t_a \neq t_b$, we define this as a time deviation.
\end{definition}

We employ a hierarchical branching scheme and branch initially on station deviations and, if necessary, subsequently on time deviations.
We introduce branching constraints by modifying the graphs' arc sets $\mathcal{A}_v$ to either enforce or forbid visiting charging stations entirely for station deviations and specific charging vertices for time deviations. While branching, we maintain a global column pool, which we update with all new columns found. We use the subset of columns that are valid for all activated branching constraints as a local set of columns in each branching node. We use a depth-first search strategy to find a feasible solution and branch on a single vehicle at a time. 
To search for an integer solution in the existing column pool, we periodically solve a branching node under restored integrality constraints using a commercial integer solver with a brief time limit of one minute to utilize its optimized branch-and-bound routines and heuristics.

\subsection{Robust Solution Framework}\label{sec:robust_master}\label{sec:robust_master:relaxed}

In order to consider uncertain \revnew{customer schedules}, we embed our deterministic algorithm into a robust solution framework. We consider an adjustable robust problem setting \citep[cf.][]{ben-tal_adjustable_2004} that formally reads
\begin{equation}\label{eq:general_adjustable_robust}
\begin{array}{ll@{}ll}
\text{minimize}  &   \boldsymbol{c}^T (\boldsymbol{u}^T, \boldsymbol{v}^T(\boldsymbol{\zeta}))  &  &\\
\text{subject to}&  \boldsymbol{A}(\boldsymbol\zeta) (\boldsymbol{u}^T ,\boldsymbol{v}^T(\boldsymbol\zeta))   \leq \boldsymbol{b} &\ \forall\boldsymbol{\zeta}\in\mathcal{U},\\
&                 \boldsymbol{c} \in \mathbb{R}^n, \boldsymbol{A} \in \mathbb{R}^{m \times n}, \boldsymbol{b} \in \mathbb{R}^{m}&\\
                 
\end{array}
\end{equation}
 with a vector of objective function coefficients $\boldsymbol{c}$, a vector of adjustable decision variables 
$\boldsymbol{v}(\boldsymbol\zeta)$ and a vector of unadjustable decision variables $\boldsymbol{u}$, a right-hand side vector $\boldsymbol{b}$,
a left-hand side matrix   $\boldsymbol{A}$, and an uncertainty set $\mathcal{U}$. A decision on  $\boldsymbol{u}$ must be taken before the 
realization $\zeta$ of the uncertainty set $\mathcal{U}$ is known, while the adjustable decision variables $\boldsymbol{v}(\boldsymbol\zeta)$ can be modified afterward.

In our setting, the unadjustable decision is to select the charging stations, while operational routing and charging decisions can be adjusted once the uncertainty realization is known. By design, our decomposition of the deterministic planning problem separates adjustable, and unadjustable decisions: unadjustable network design decisions remain in the \gls{abk:cpp}, while adjustable operational decisions remain in the second stage operational problem. 

\revnew{The uncertainty results from varying vehicles' customer schedules and accordingly yields varying driving patterns.}
Since such variations cannot be modelled with standard ellipsoidal-uncertainty sets, we approximate $\mathcal{U}$ with a primitive uncertainty set $z\in\mathcal{Z}$ which contains specific scenarios that we derive from real-world data. We detail the generation of these scenarios in Section~\ref{sec:computational_study:customer_demand_data} and  verify the quality of the approximation through a-posteriori validation on a separate set of out-of-sample uncertainty realizations in Section~\ref{subsec:validation_of_uncertainty_approx}.

In the following, we first introduce a \gls{abk:FSA} that considers the entire scenario set~$\mathcal{Z}$. We then introduce modified approaches that consider only subsets of scenarios under relaxed feasibility requirements.

The \gls{abk:FSA} outlined in Algorithm~\ref{alg:fsa}  searches for a solution of the \gls{abk:cpp} which is feasible for all scenarios.
Specifically, we solve the \gls{abk:cpp} (l.\ref{lst:fsa:bmp}) and subsequently evaluate a separate second stage operational problem for every scenario $z\in\mathcal{Z}$ (l.\ref{lst:fsa:oracle}). If the operational problem of any scenario is infeasible, we create a cover cut and continue with the next \gls{abk:cpp} iteration (l.\ref{lst:fsa:add_cut}).
We modify the cover strengthening heuristic such that it requires operational feasibility for all scenarios before reducing the cover size  (l.\ref{lst:fsa:lift}).

\begin{algorithm}
\caption{\gls{abk:FSA}}\label{alg:fsa}
\begin{algorithmic}[1]
\linespread{1}\selectfont
\State Initialize \gls{abk:cpp}\label{lst:fsa:ccp}
    \Repeat
        \State $\overline{x}_s \gets $ Solve \gls{abk:cpp}\label{lst:fsa:bmp}
        \ForAll{$z \in \mathcal{Z}$}
            \If {operational problem for scenario $z$  is not feasible given $\overline{x}_s$}\label{lst:fsa:oracle}
                \State $\mathcal{I} \gets \{ s \in \mathcal{N}^s\ |\ \overline{x}_s = 0 \}$\label{lst:fsa:create_cover}\;
                \State $\mathcal{O} \gets $ \Call{Strengthen Cover$_\mathcal{Z}$}{$\mathcal{I}$}  \Comment{mod. Algorithm~\ref{alg:lifting}}\label{lst:fsa:lift} \;
                \ForAll{$\mathcal{P} \in \mathcal{O}$}
                \State Add constraint $\sum_{s \in \mathcal{P}} x_s \geq 1$ to \gls{abk:cpp}\label{lst:fsa:add_cut}\;
                \EndFor
                \State \textbf{break}
            \EndIf
        \EndFor
    \Until{All scenario operational problems confirm feasibility}\label{lst:fsa:exit}
\end{algorithmic}
\end{algorithm}

The \gls{abk:FSA} enforces full feasibility for all considered scenarios.
In practice, a fleet operator might accept a small amount of infeasibility to obtain lower costs.
We can quantify the feasibility of a solution using two different perspectives by assessing either the vehicle feasibility or the scenario feasibility.
\begin{definition}[vehicle feasibility]
    Given a charging station configuration $\mathcal{S}$ and the set of vehicles $\mathcal{V}_z$ in scenario $z$, we define the set of feasible vehicles $\mathcal{V}^f_z$ as the maximum size subset of $\mathcal{V}_z$ such that each vehicle  in $\mathcal{V}^f_z$ is operationally feasible for the installed charging stations $\mathcal{S}$. We then define mean vehicle feasibility  as \begin{equation}\overline{V}_\text{feas} = \frac{1}{|\mathcal{Z}|} \sum_{z \in \mathcal{Z}}\frac{|\mathcal{V}_z^f|}{|\mathcal{V}_z|}\end{equation} and the minimum (worst case) vehicle feasibility as 
    \begin{equation} 
            V_\text{feas}^{\min} = \min \left\{ \frac{|\mathcal{V}_z^f|}{|\mathcal{V}_z|}\ |\ z \in \mathcal{Z} \right\}.
    \end{equation}
\end{definition}
\begin{definition}[scenario feasibility]
Given a charging station configuration $\mathcal{S}$ and a set of scenarios $\mathcal{Z}$, we define the scenario feasibility $Z_\text{feas}$ of  $\mathcal{S}$ as the percentage of scenarios $z \in \mathcal{Z}$ in which all vehicles are feasible, i.e., \begin{equation}
    Z_\text{feas} = \frac{1}{|\mathcal{Z}|} \sum_{z \in \mathcal{Z}}
    \begin{cases} 
    1, & \text{if}\ |\mathcal{V}^f_z| = |\mathcal{V}_z| \\
    0, & \text{otherwise}
    \end{cases}\end{equation}
\end{definition}

We developed a parameterized solution approach for each feasibility type, which allows to obtain a solution for a certain feasibility level $\alpha$.

\paragraph{\gls{abk:aSA}}
This approach is motivated by the concept of $\alpha$-reliability proposed by \cite{daskin_-reliable_1997}, whose model self-selects a ``reliability set'' of scenarios such that their cumulative likelihood of occurrence is above a given $\alpha$.
We design our \gls{abk:aSA} to find an optimal solution  feasible for at least $\alpha$ of all scenarios, as all our scenarios have an equal probability of occurrence. 
In order to do so, we modify Line~\ref{lst:fsa:exit} of Algorithm~\ref{alg:fsa} to accept a solution as soon as it is feasible for at least $\alpha$ percent of the scenarios.
In contrast to the \gls{abk:FSA}, the \gls{abk:aSA} cannot terminate early once one scenario is detected to be infeasible but must carry on until $1-\alpha$ infeasible scenarios are found.

\paragraph{\gls{abk:aVA}}

In the \gls{abk:aVA}, we seek for a solution in which at least $\alpha$ percent of the vehicles are fully operational in each scenario.
Here, instead of applying the feasibility relaxation on the aggregated scenario set level, we relax the feasibility requirements for accepting an individual scenario. Accordingly, we modify the \gls{abk:FSA} to accept a scenario as operationally feasible if we can obtain a vehicle feasibility of at least $\alpha$. To calculate this vehicle feasibility, we need to compute the number of feasible vehicles.
Within the column generation master problem, we already obtain a separation of vehicles into feasible and infeasible vehicles through the selection of dummy routes $y_{v0}$. However, so far we can exit the procedure early once the first infeasible vehicle is detected.
In order to bound the vehicle feasibility to some $\alpha$-level, we need to determine the smallest subset of infeasible vehicles. This takes considerable computation time, especially during the cover strengthening heuristic, which repeatedly evaluates the operational problem.
Hence, we base our \gls{abk:aVA} approach on the scenario-based adversarial sampling method introduced by~\cite{bienstock_computing_2008}. 
The latter has initially been designed to tackle an optimization problem where the number of scenarios is too large to be incorporated into a single model. 
The set of all scenarios is divided into disjoint subsets: The optimization set and the feasibility set. The idea is to derive optimality for each scenario in the optimization set and feasibility for each scenario in the feasibility set. 
\citeauthor{bienstock_computing_2008} select a small subset of scenarios to seed an optimization set, which they use to obtain an optimal solution.  They then evaluate if this solution is at least feasible for all remaining scenarios. Afterward, they add any infeasible scenarios to the optimization set and repeat the procedure until the solution is optimal for the optimization set and at least feasible for the remaining scenarios.
We adapt this approach to our problem setting by introducing a feasibility set of scenarios for which we guarantee full feasibility $\Omega$, and a separate
relaxed feasibility set $\Phi = \mathcal{S} \setminus \Omega$ for which we guarantee at least $\alpha$-level feasibility.
This separation allows us to use the faster \gls{abk:FSA} method on the full feasibility set $\Omega$, and only requires to compute 
the exact number of feasible vehicles when we evaluate the $\alpha$-level acceptance function for the scenarios in the relaxed feasibility set $\Phi$.

Algorithm~\ref{alg:bas} shows our \gls{abk:aVA} approach.
First, we initialize the feasibility set $\Omega$ with a single seed scenario (l.2) and apply the \gls{abk:FSA} to $\Omega$ to obtain $\overline{x}_s$  (l.5--11). 
We then compute the vehicle feasibility of each scenario $z \in \Phi$ (l.13) by solving the operational problems for the respective $\overline{x}_s$ and counting vehicles for which a non-dummy route was selected. 
If each scenario of the relaxed feasibility set $\Phi$ is within the $\alpha$-level of feasibility, we can terminate the algorithm (l.16). Otherwise, we select the scenario $z^*$ with the smallest $\alpha$-level violation as a candidate for the full feasibility set $\Omega$ (l.13). Adding the scenario $z^*$ to $\Omega$, and thus including it in the \gls{abk:FSA}, would result in 100\% feasibility for all vehicles in $z^*$.  However, especially at lower $\alpha$-levels, this would be an over-fulfillment of vehicle feasibility for scenario $z^*$. To reduce this potentially costly over-fulfillment, we derive a partial scenario $z'$ which contains a subset of the vehicles contained in scenario candidate $z^*$ (l.14). Here, we include only those vehicles required to cover the difference between the evaluated vehicle feasibility and the desired $\alpha$-level. For example, if the chosen $\alpha$-level allows 50 infeasible vehicles, but the current number in scenario $z^*$ is 60, we create a scenario $z'$, which consists of 10 randomly selected infeasible vehicles of $z^*$. We then add $z'$ to $\Omega$ (l.15). 
We repeat this process until the vehicle feasibility for each scenario in the relaxed feasibility set $\Phi$ is larger or equal to $\alpha$ (l.20). In contrast to~\cite{bienstock_computing_2008}, we add only one scenario to the feasibility set $\Omega$ in each iteration. 
The final solution is entirely scenario feasible over the optimization set and at least  $\alpha$-level vehicle feasible for all other scenarios. 

The seed scenario has a significant impact on the behavior of the \gls{abk:aVA} algorithm. Therefore, we evaluate two selection strategies for the seed scenario based on the cost of the solutions derived by independently solving each scenario as a deterministic problem. For the first selection strategy, we pick the scenario with the lowest cost, while for the second selection strategy, we select the scenario with the median cost.  We use the cuts generated for the seed scenario to warm-start the \gls{abk:aVA}. In the remainder of this paper, we refer to a specific version of the \gls{abk:aVA} by stating its $\alpha$-level and suffixing the letter \emph{L} for the lowest cost seed, or the letter \emph{M} for the median cost seed, e.g., 95-VA-L is the notation for the vehicle feasibility approach with a feasibility level of $\alpha=95\%$ and the lowest cost scenario as seed.

\begin{algorithm}
\caption{\glsfirst{abk:aVA}}\label{alg:bas}
\begin{algorithmic}[1]
\linespread{1}\selectfont
\State Initialize \gls{abk:cpp}\label{lst:bas:init}
\State Add seed scenario to $\Omega$
\Repeat
    \Repeat
        \State $\overline{x}_s \gets $ Solve \gls{abk:cpp}\label{lst:bas:bmp}
        \ForAll{$\omega \in \Omega$}
            \If {operational problem for scenario $\omega$  is not feasible given $\overline{x}_s$}\label{lst:bas:oracle}

            \State $\mathcal{I} \gets \{ s \in \mathcal{N}^s\ |\ \overline{x}_s = 0 \}$\label{lst:bas:create_cover}\;
            \State $\mathcal{O} \gets $ \Call{Strengthen Cover$_\Omega$}{$\mathcal{I}$}  \Comment{mod. Algorithm~\ref{alg:lifting}}\label{lst:bas:lift} \;
            \ForAll{$\mathcal{P} \in \mathcal{O}$}
            \State Add constraint $\sum_{s \in \mathcal{P}} x_s \geq 1$ to \gls{abk:cpp}\label{lst:bas:add_cut}\;
            \EndFor
                \State \textbf{break}
            \EndIf
        \EndFor
    \Until{$\overline{x}_s$ feasible for all scenarios in $\Omega$}
        \State $z^* \gets$ infeasible scenario $ z \in \Phi$ with minimal violation
        \State $z' \gets$ derive partial scenario from $z^*$
        \State Add $z'$ to $\Omega$\label{lst:bas:add}
\Until{all scenarios in $\Phi$ accepted with at least $\alpha$-level vehicle feasibility}
\end{algorithmic}
\end{algorithm}

\revnew{
\subsection{Selection and Sizing of Charging Stations}\label{sec:variable_sizing}
Motivated by discussions with practitioners, our basic algorithmic framework is limited to a fixed parameter $n_s$ for the number of charging points. In practice, such an approach is favorable for planning simplicity, i.e., requiring only a single technical specification and easing the analysis of potential station locations, but also for equity constraints, i.e., offering the same amount of chargers throughout a city when building a new charging station. Lastly, ordering a larger number of homogeneous charging stations yields cost advantages in terms of discounts for the municipality. Nevertheless, making a sizing decision within the optimization model, i.e., turning $n_s$ into a decision variable, remains interesting from a scientific perspective. Accordingly, we introduce a variant of our algorithm that allows for a variable number of charging points. We refer to this variant as ``variable-size'' and present an adapted version of the cutting plane master problem and of the cover strengthening heuristic, respectively.

As the range of the number of charging points is restricted by physical constraints, we discretize it into
an ordered set of charging point levels $\mathcal{L}$. We introduce variable $x_{sl}$, which represents the opening of a charging station $s$ at level $l$ with $n_{l}$ charging points at a cost of $o_{ls}$. Levels in $\mathcal{L}$ are ordered by the number of charging points. Our cuts $\forall j \in \mathcal{J}$ then consist of a tuple $(s,l)$ denoting the site $s$ and level $l$ included in the cut. We account for the hierarchical relationship between levels (e.g., if a site of the cover is currently opened at level 2, opening it at any larger level could potentially resolve the infeasibility) during the cut generation as follows. 
If the cover $\mathcal{I}_j$ from which the cut is constructed includes site $s_1$ at a charging point level $l_1$, we include
tuples for level $l_1$ and all larger levels, $(s_1, l) \in \{ l \in \mathcal{L}\ |\ l \geq l_1\}$
in the cut.

Our cutting plane problem for variable sizes of charging stations then reads:

\begin{mini!}
    {}{\sum_{l \in \mathcal{L}} \sum_{s \in \mathcal{N}^s}c_{sl} \cdot x_{sl}\label{eq:benders_variable_master:obj}}{\label{eq:benders_variable_master}}{}
    \addConstraint{\sum_{ (s,l) \in \mathcal{I}_j} x_{sl} }{\geq 1}{\quad \forall j \in \mathcal{J}\label{eq:benders_variable_master:cuts}}
    \addConstraint{ \sum_{l \in \mathcal{L}} x_{sl} }{\leq 1}{\quad \forall s \in \mathcal{N}^s\label{eq:benders_variable_master:convex}}
    \addConstraint{ x_{sl}                          }{\in \{0,1\},}{\quad \forall s \in \mathcal{N}^s, l \in \mathcal{L}.\label{eq:benders_variable_master:domain}}
\end{mini!}

We also adapt the strengthening heuristic to consider charge point levels.
Algorithm~\ref{alg:lifting_specific} shows the modified strengthening heuristic. We take an array $p$ containing the charging point level for each site $s$ as input. In each iteration, we pick a subset of charging stations that are not at their maximum charging point level and increase their level by one. We have successfully strengthened the cut if the resulting charging station configuration remains operationally infeasible.

\begin{algorithm}
    \caption{Variable-Size Cut Strengthening Heuristic}\label{alg:lifting_specific}
    \linespread{1}\selectfont
    \begin{algorithmic}[1]
        \Function{Strengthen Cut}{$p$} \label{lst:line:cover_lift}
        \State $\mathcal{O} \gets \{\Call{DeriveCut}{p}\}$
        \While{time limit not exceeded}\label{lst:line:time_limit}
          \State $\tilde{p} \gets p$
           \While{$ \exists \{ l < \max(\mathcal{L}) | l \in p \}  $} 
   
               \State $\tilde{p} \gets$ increase levels of random subset of $\tilde{p}$\label{lst:line:random_subset}
               \If{ levels $\tilde{p}$ are not operationally feasible}
                   \State $\mathcal{O} \gets \mathcal{O} \cup \{\Call{DeriveCut}{\tilde{p}\}}$\label{lst:line:add_to_O}
               \Else
                   \State break\label{lst:line:break_on_feasible}
               \EndIf
           \EndWhile
          \EndWhile
          \State Prune $\mathcal{O}$ to retain only non-dominated cuts\;
          \State \Return $\mathcal{O}$
        \EndFunction
         \end{algorithmic}
    \end{algorithm}

Introducing variables for all charging point levels increases the total number of decision variables in the strategic problem by $(|\mathcal{L}| - 1) \cdot |\mathcal{N}^s|$. We partially mitigate this complexity increase by employing a two-step solution method. First, we note that all cuts generated in the original fixed-size model must also be true in the variable-size model, i.e., if opening a charging station with $n_s$ charging points does not resolve an infeasibility, neither would an even smaller charging station size. We can therefore solve the static-size model using the largest size level in $\mathcal{L}$ and then use the generated cuts to warm-start the variable-size model. As a second improvement, we do not add all columns with negative reduced cost generated by the pricing problem to the column pool to limit its size. If we generate more than 50 columns for a vehicle, we apply a $k$-metroid based sampling with a similarity measure as a distance function to pick 50 diverse columns. We measure the similarity of two columns by counting the periods where the two columns use the same charging site. Third, we modified the operational branch-and-price problem by including a more aggressive diving heuristic. We sort the vehicles by the number of fractional columns selected. If the number of columns exceeds a threshold, we add a branch that fixes the first column for the first vehicle (i.e., the one with the largest amount of fractional columns). We then dive on these branches until we reach feasibility. More aggressive diving heuristics, such as the ``Diving for Feasibility'' heuristic using Limited Discrepancy Search proposed in \cite{sadykov_primal_2019}, did not result in better solution times for our problem. We note that the presented modifications did not result in consistent improvements for the static-size case and, thus, were only used for the variable-size case.



}


\section{Experimental Setup}\label{sec:experimental_design}

In order to illustrate the advantages of our proposed algorithms, we created several instances of varying sizes from real-word data.
Further, we introduced a special set of battery and charger type specific instances to generate insight into how technological change may impact charging network design.

\subsection{Real-World Data}

We use GPS data from a taxi association in the metropolitan area of Munich, Germany, collected by \cite{jager_analyzing_2016} for our study.
The data set contains all customer trips undertaken by a fleet of 500 taxis powered by combustion engines during 2015.  
The number of active taxis per day ranges between 300 and 500. We complemented this data with information about the location of existing taxi stands, and calculated distances and travel times from customer locations to charging stations \revnew{using the Open Source Routing Machine project \citep[cf.][]{luxen-vetter-2011} and \citeauthor{openstreetmap_openstreetmap_2016} road network data.}

\begin{figure}[h!] \FIGURE {\includegraphics[]{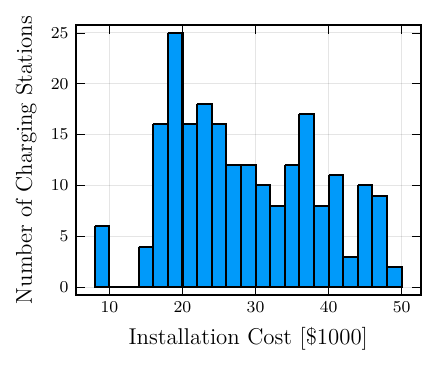}}
{Distribution of the \revnew{opening} cost for 215 charging stations\label{fig:cost_function}} {} \end{figure}

\paragraph{Locations \revnew{and Sizing} of Charging Stations}
We identified charging station locations based on the existing 215 taxi stands of the city of Munich. 
\revnew{
    We selected $n_s=4$ charging points for every charging location, as larger charging station configurations have cost advantages in construction and maintenance. The taxi stands in the city have an average number of 7.2 parking spaces, leaving enough space for charging operations of that size. The selection of $n_s=4$ also aligns with the procedures of the largest local charging station operator in Munich, which primarily constructs stations with four charging points. For the opening costs, we assume that locations closer to the city center have higher costs. Using the normalized distance to the city center,  $d$, $0 \leq d \leq 1$, we derive the cost function: }
\begin{equation}\label{eq:site_cost}
    \revnew{o}(d) = 30 \cdot e^{-10d} - 10d + 20,
\end{equation}
\revnew{where the parameters have been set based} on the empirical study of \cite{smith_costs_2015}. The minimum costs are $\revnew{o}(1)=10$ at the city limits, and the maximum costs are $\revnew{o}(0)=50$ in the city center, with an exponential increase of costs when approaching the city center. Figure~\ref{fig:cost_function} shows the distribution of the \revnew{opening} costs for the 215 stations.

\revnew{
    In the variable-size approach, we chose between two or four charging points since charging stations available on the market usually provide an even number of charging points. Dedicating more than four spaces at a taxi stand for charging (more than half the average number of parking spaces) would hinder normal taxi operations.  Following \cite{nicholas_estimating_2019}, we set the cost of a two-charging-point station to be 67\% of the cost of a four-charging-point station.
}

\paragraph{\revnew{Vehicles' Customer Schedules}}\label{sec:computational_study:customer_demand_data}

In our analysis we consider 24/7 operations which requires to charge taxis in between customer trips during operations. Additionally, we assess a depot charging setting in Section~\ref{sec:comp:tech_change}. We use anonymized data and thus cannot identify taxis over multiple shifts. Therefore, we model the customer schedule of each taxi shift as an individual vehicle. To model handover between shifts, we set for each vehicle $v$ the initial state of charge $q^\text{b}_v$ and the state of charge required at the end of the shift $q^\text{e}_v$ to 50\% of its maximum battery capacity $q^{\max}$.
We removed taxis with shifts below 5 hours, as these taxis are not operating in a multi-shift system and are therefore capable of using (depot) charging stations not considered in our study.
Table~\ref{tbl:instance_stats} provides key characteristics of the taxis for a 24-hour time period. The shift duration median is 8 hours, reflecting the local regulations limiting drivers' duty times. 

\begin{table}[h!]
    \TABLE
    { Characteristics for taxi shifts per day.\label{tbl:instance_stats}}
    {
    \begin{tabular}{llllll}
    \toprule
    {} &                    Mean &                     Std &              Min &             Median &              Max \\
    \midrule
    Shift duration [h:m:s]    &  07:41:44 &  00:56:31 &  04:57:46 &  08:00:15 &  08:54:27 \\
    Driven distance [km]    &                  123.06 &                 46.97 &            30.22 &           118.43 &           250.28 \\
    Number of customers &                 10.92 &                 2.73 &                4 &               11 &               19 \\
    \bottomrule
    \end{tabular}
    }
    {}
    \end{table}

\begin{figure}[h!] \centering
        \begin{subfigure}{0.4\textwidth}
        \includegraphics[width=\textwidth]{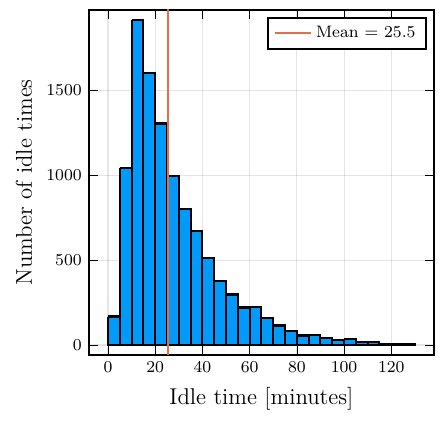}
        \caption{Idle Time between Customer Trips}
        \end{subfigure}
        \begin{subfigure}{0.4\textwidth}
        \includegraphics[width=\textwidth]{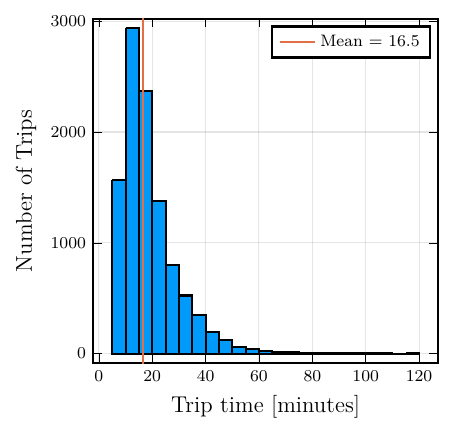}
        \caption{Duration of Customer Trips}
        \end{subfigure}
        
        \caption{Distribution of durations within taxi shift for the first scenario of the instance with $|V|=1,000, |N^s| = 60$}\label{fig:shift_timings}
    \end{figure}

    While 8 hour shifts would leave some room for charging at a depot, we assume that taxis are handed over between shifts and thus operate 24/7. Thus, charging has to be undertaken during operations between customer trips. 
    \revnew{ 
        The static customer sequence assigned to each vehicle specifies the time windows under which charging might occur. Our data shows that the average idle time between two customers is 25 minutes, which is sufficient to enable flexible charging operations. Figure~\ref{fig:shift_timings} shows the empirical distributions of  the idle time and the duration of customer trips. The idle time is largely available for charging as the local taxi ordinance requires taxis to wait at stands and prohibits them from actively soliciting customers or repeatedly circling within an area.
    }

\paragraph{Vehicle, Battery, and Charging Parameters}
We assume a fleet of e-NV200-Taxis, an electric taxi from Nissan based on the popular Nissan Leaf platform with a 40 kWh battery and 250 km range.
We assume a charging speed of 50 kW and a constant-current constant-voltage charging process. 
We derive a non-linear charging function from the methodology and data presented by \cite{strasser_synthetische_2013}.
Under these parameters, 95\% of all observed taxi shifts are operationally feasible when making all charging stations available without capacity restrictions.

\subsection{Instances}

We created instances that are parameterized by the number of vehicles and the number of potential charging stations. We vary the number of vehicles between 100 and 1,500, and the number of charging stations between 30 and 70.  We used the instance with 1,000 vehicles and 60 charging stations to analyze the robust algorithms, as this instance represents a reasonable electrification target for the city of Munich.
\revnew{Each instance consists of 28 different scenarios, modeling a month of customer demand where each scenario represents the vehicles' customer schedules from 24-hours of operations}. We created scenarios for an instance sized at a specific number of vehicles by applying proportional stratified sampling.
We used all taxi shifts of a 6-month time window and grouped them by the day-of-the-week of their start time. For each group, we stratified the shifts by their start time and proportionally sampled from every stratum.
We sample four scenarios for each day of the week to form the instances' scenario set.
Where shifts extend into the following day, we mapped the charging demand into the day under consideration.   
In order to vary the number of charging station locations for different instance sizes, we apply $k$-Means clustering. With $k$ as the number of desired potential charging station locations, we obtained $k$ cluster centroids and for each centroid, choose the charging station closest to it. This method distributes the charging station locations evenly spaced throughout the city for any given instance size.
As we vary both the number of available charging station locations and the number of vehicles independently, we potentially create instances 
with either a large mismatch between demand and supply or instances where some vehicles have no reachable charging stations. We remove the vehicles that would induce infeasibility and sample replacement vehicles for those instances.

We created a special set of instances for studying the impact of technological parameters.  Varying the technological parameters influences vehicle feasibility significantly, leading to many infeasible scenarios, especially in case of small battery-sizes. For these settings, we might be unable to obtain solutions using the FSA as it requires all scenarios to be potentially feasible. We resolved this problem by generating parameter specific instances:
We first generated a master instance with 500 vehicles and 60 charging stations using the original sampling method but with a charging speed and battery size twice as high as the baseline configuration when resampling vehicles. This master set thus contains those vehicles that would be operationally feasible assuming considerable advances in technology and the opening of all available charging stations. From this master set, we then derive parameter-specific instances for every charging and battery parameter variation by removing vehicles that are not operationally feasible under the selected parameters, even if we open all charging stations. Our parameter-specific instances thus contain the subset of vehicles that could complete their shifts when considering the given battery and charging parameters.
When evaluating a solution for vehicle feasibility, we count all vehicles removed for the derived instances as infeasible and do not evaluate the scenario feasibility in this setting.


\section{Results}\label{sec:results}
In this section, we show the results for our experiments. We first focus on the performance evaluation of the deterministic algorithm using the benchmark instances, and continue by evaluating the proposed robust algorithms. Finally, we focus on a detailed solution analysis and explore the variation of charging technological parameters.

All algorithms were implemented in Rust using Gurobi 9.0.3 as linear programming solver and were run on a Linux 5.17 operating system using a single core of an Intel i7-6700 processor with 3.40GHz. 

\subsection{Computational Performance}\label{sec:performance_evaluation}

In order to assess our general algorithmic approach presented for the deterministic problem in Section~\ref{sec:solution_method}, we solve each scenario of an instance independently. We observe the computation time, costs, and the number of opened sites and report aggregated results for each instance. We vary the number of vehicles $|V|$ between 100 and 1,500, and the number of charging stations $|\mathcal{N}^s|$ between 30 and 70.
Table~\ref{tbl:results_deterministic} states the average computation time, costs, and the number of opened charging stations. As can be seen, computation time, cost, and the number of opened charging stations increase with the size of the fleet and the number of available stations, respectively. 
\begin{table}\small
    \TABLE
    {Solutions for the deterministic model for varying instances stated as distributions over 28 scenario realizations.\label{tbl:results_deterministic}}
    { \begin{tabular}[t]{rrrrr}
\toprule
 \multicolumn{2}{c}{Instance}  & & & \tabularnewline
\cmidrule{1-2}
$|V|$ & $|\mathcal{N}^s|$  & $t$ & $o$ & $|S^*|$\tabularnewline
\midrule
100&30&24.4&264.0&7.7 \tabularnewline
200&30&34.0&363.4&10.5 \tabularnewline
300&30&41.5&434.7&12.6 \tabularnewline
400&30&47.7&497.9&14.4 \tabularnewline
500&30&54.4&539.1&15.6 \tabularnewline
600&30&59.4&583.2&16.9 \tabularnewline
700&30&63.9&621.4&18.0 \tabularnewline
800&30&79.5&655.4&18.9 \tabularnewline
900&30&87.4&686.4&19.9 \tabularnewline
1000&30&103.2&715.8&20.8 \tabularnewline
1100&30&178.9&725.4&21.0 \tabularnewline
1200&30&316.9&746.2&21.6 \tabularnewline
1300&30&532.6&766.5&22.0 \tabularnewline
1400&30&1,066.6&780.4&22.5 \tabularnewline
1500&30&2,408.3&789.0&22.7 \tabularnewline
\bottomrule
\end{tabular}
\begin{tabular}[t]{rrrrr}
\toprule
 \multicolumn{2}{c}{Instance} &  & & \tabularnewline
\cmidrule{1-2}
$|V|$ & $|\mathcal{N}^s|$ &   $t$ & $o$ & $|S^*|$\tabularnewline
\midrule
100&50&100.1&289.6&8.4 \tabularnewline
200&50&109.9&448.3&12.7 \tabularnewline
300&50&138.4&554.1&15.6 \tabularnewline
400&50&123.9&633.9&17.6 \tabularnewline
500&50&146.4&712.7&19.9 \tabularnewline
600&50&165.4&754.6&21.1 \tabularnewline
700&50&171.5&815.7&22.7 \tabularnewline
800&50&204.2&857.7&23.9 \tabularnewline
900&50&235.0&903.7&25.4 \tabularnewline
1000&50&279.6&943.8&26.5 \tabularnewline
1100&50&317.1&976.3&27.4 \tabularnewline
1200&50&432.9&1,008.6&28.5 \tabularnewline
1300&50&649.1&1,046.4&29.5 \tabularnewline
1400&50&1,012.8&1,078.0&30.5 \tabularnewline
1500&50&1,227.0&1,095.5&31.1 \tabularnewline
\bottomrule
\end{tabular}
\begin{tabular}[t]{rrrrr}
\toprule
 \multicolumn{2}{c}{Instance} &  & & \tabularnewline
\cmidrule{1-2}
$|V|$ & $|\mathcal{N}^s|$ &   $t$ & $o$ & $|S^*|$\tabularnewline
\midrule
100&70&586.4&323.2&9.4 \tabularnewline
200&70&938.8&461.9&13.1 \tabularnewline
300&70&1,339.1&567.0&16.3 \tabularnewline
400&70&833.5&654.0&18.7 \tabularnewline
500&70&930.8&728.5&21.1 \tabularnewline
600&70&837.0&792.0&22.8 \tabularnewline
700&70&1,133.8&851.5&24.4 \tabularnewline
800&70&1,293.4&901.2&25.9 \tabularnewline
900&70&1,105.6&962.3&27.6 \tabularnewline
1000&70&1,483.6&998.1&28.8 \tabularnewline
1100&70&1,594.1&1,039.3&29.9 \tabularnewline
1200&70&2,436.0&1,077.6&31.1 \tabularnewline
1300&70&3,942.4&1,109.5&31.9 \tabularnewline
1400&70&4,513.3&1,143.1&33.0 \tabularnewline
1500&70&5,847.9&1,172.8&33.9 \tabularnewline
\bottomrule
\end{tabular}
}
    {Abbreviation: $|V|$ - avg. number of vehicles, $|\mathcal{N}^s|$ - number of potential charging stations, $t$ - computation time in seconds, $o$ - cost of solution , $|S^*|$ avg. number of opened charging stations }
    \end{table}

\begin{figure}[h!]
    \FIGURE
    {\includegraphics[]{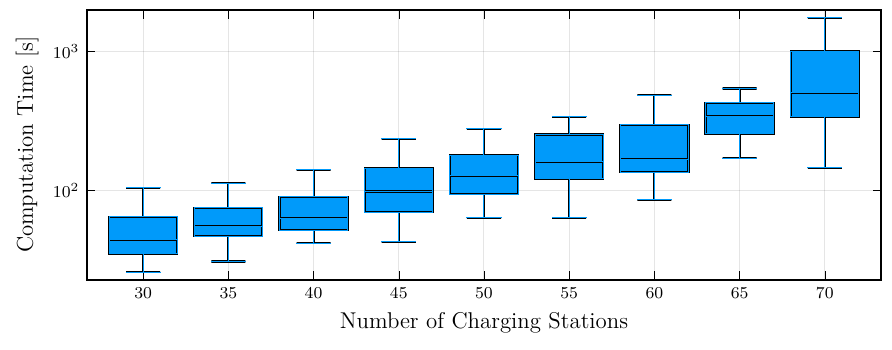}}
    {Computation time (log scale) of the algorithm for instances by the number of charging stations for $|\mathcal{V}| = 500$. \label{fig:boxresults:duration_byvehicle_500}}
    {}
    \end{figure}

\begin{figure}[h!]
        \FIGURE
        {\includegraphics[]{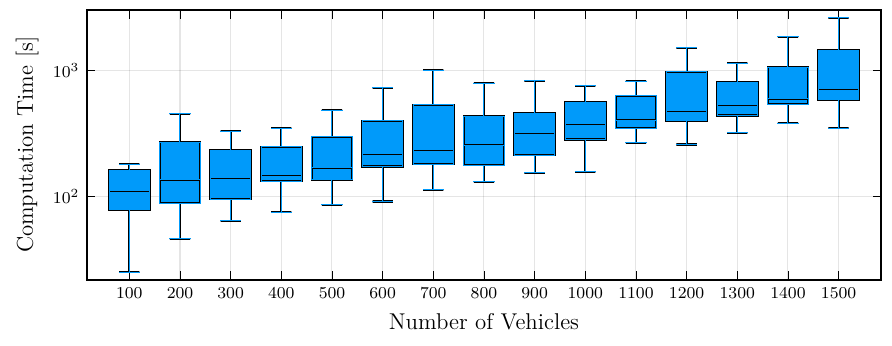}}
        {Computation time (log scale) of the algorithm 
            for instances by number of vehicles for $|\mathcal{N}^s| = 60$. \label{fig:boxresults:duration_bysize_60}}
            {}
    \end{figure}

\begin{figure}[h!]
        \FIGURE
        {\includegraphics[]{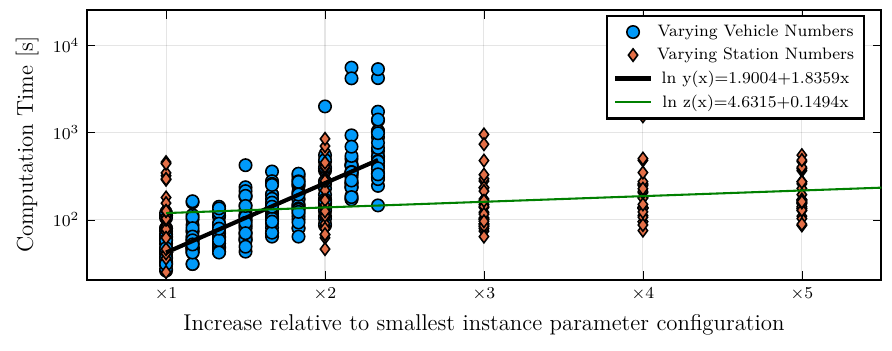}}
        {
             Individual scenario computation times for increasing number of charging stations under constant number of vehicles and vice versa. The x-axis shows the normalized increase relative to each respective smallest configuration.\label{fig:boxresults:relative_increase}}
            {}
    \end{figure}

To isolate the influence of each of the instance parameters, we use a ceteris-paribus sensitivity study. We show box-whiskers plots of the computation time distributions when varying the number of potential charging stations under fixed fleet size in Figure~\ref{fig:boxresults:duration_byvehicle_500}, and increasing vehicle sizes under a fixed number of charging stations in Figure~\ref{fig:boxresults:duration_bysize_60}. 
We observe a strong exponential increase in the median computational time and variance when increasing the number of potential charging stations between 30 and 70 under a fixed fleet size of 500. While the median computation time is 54 seconds for 30 charging stations, it increases to 15 minutes for 70 charging stations. When varying the number of vehicles between 100 and 1,500, we still observe a roughly exponential trend but with a lower rate of increase than in Figure~\ref{fig:boxresults:duration_byvehicle_500}. 
Note, however, that the parameter range for each type is different as they were selected to contain reasonable real-world configurations.  The largest number of vehicles is about 15 times larger than the smallest, while the largest number of charging stations is only about 2.3 times as large as the smallest.
We show a normalized comparison between both methods using a logarithmic transformation in Figure~\ref{fig:boxresults:relative_increase}. The solution times of the individual scenarios are depicted in a scatter plot, while the two lines show the result of applying a log-level regression.  The x-axis represents the relative change compared to the smallest value of each parameter range. The regression coefficients show that a 100\% increase in the number of charging stations leads to a 527 \% increase in computation time, while a 100\% increase in the number of vehicles leads only to a 16\% increase in computation time. The number of potential charging stations directly affects the combinatorial cutting plane master problem, and thus results in a substantial increase in computation time. In contrast, the number of vehicles only affects the operational second stage problem. In our decomposition, each vehicle has its own pricing problem, and the interaction between vehicles is limited to the coupling constraints, such that the increase in computation time per vehicle remains limited.


\subsection{Robust Charging Network Design}\label{sec:robust_evaluation}

After evaluating the deterministic algorithm's performance characteristics, we evaluate the robust solution methods to show their computational performance and the behavior of the relaxed feasibility methods for varying $\alpha$-levels. We use the instance with 1,000 vehicles and 60 charging stations to evaluate the robust algorithms. We evaluate the robust solution methods in terms of objective function value, computation time, and feasibility. To evaluate the robustness of a solution, we use the scenario and vehicle feasibility as introduced in Section~\ref{sec:robust_master}. We compare the different approaches using different seed scenarios and $\alpha$-levels between 100\% and 50\%.
We introduce the \gls{abk:ISA} as a benchmark. Here, we solve each scenario independently for each instance, rank the scenario solutions by cost, and pick the solution with the median (\gls{abk:ISA}-M) or lowest (\gls{abk:ISA}-L) cost as the robust solution. The scenarios selected for the \gls{abk:ISA} correspond to the seed scenario variants available in the \gls{abk:aVA}. The lowest cost scenario thus represents the lower bound for all our solution approaches, as this is the cost-minimal charging station configuration feasible in at least one scenario.
For a fair comparison, we include the time required to determine the seed scenario for the \gls{abk:ISA} and \gls{abk:aVA} in their total computation time.

\subsubsection{Comparison of the Robust Algorithm Variants}

We evaluate the robust solution methods regarding objective function value, computation time, and feasibility. To evaluate the robustness of a solution, we use the scenario and vehicle feasibility as introduced in Section~\ref{sec:robust_master}. Table~\ref{tbl:robust_results} shows an overview of the solution quality of the different methods using the \gls{abk:ISA} as a benchmark. 

\begin{table}[h!]
\small
\TABLE
{Comparison of the robust approaches for instance of size $|V|=1,000$ and $ |S|=60$.\label{tbl:robust_results}}
{%
        \begin{tabular}{l@{\hskip 15pt}rr@{\hskip 15pt}r@{\hskip 15pt}rrrr}
        \toprule
        Method &   o & $\Delta$o[\%] & t[m] & $\overline{V_{\text{feas}}}$[\%] & $V_{\text{feas}}^{\max}$[\%] & $Z_{\text{feas}}$[\%]\tabularnewline
        \midrule
        [2pt]
         \gls{abk:FSA} & 1969 & 0.0  & 50.63 &  100.0 &  100.0  & 100.0 \tabularnewline[2pt]
         \gls{abk:ISA}-L & 634 & -67.8  & 276.32 &  96.19 &  90.0  & 3.57 \tabularnewline[2pt]
         \gls{abk:ISA}-M & 965 & -50.99  & 276.32 &  98.32 &  95.8  & 3.57 \tabularnewline[2pt]
         90-SA & 1888 & -4.11  & 201.80 &  99.99 &  99.9  & 89.29 \tabularnewline
         80-SA & 1810 & -8.08  & 439.12 &  99.95 &  99.5  & 78.57 \tabularnewline
         50-SA & 1637 & -16.86  & 3349.50 &  99.86 &  99.3  & 50.0 \tabularnewline[2pt]
         99-VA-M & 1547 & -21.43  & 307.07 &  99.76 &  99.4  & 10.71 \tabularnewline
         97-VA-M & 1126 & -42.81  & 286.68 &  98.86 &  97.0  & 3.57 \tabularnewline
         95-VA-M & 965 & -50.99  & 277.42 &  98.32 &  95.8  & 3.57 \tabularnewline
         90-VA-M & 965 & -50.99  & 277.42 &  98.32 &  95.8  & 3.57 \tabularnewline[2pt]
         99-VA-L & 1518 & -22.91  & 577.97 &  99.66 &  99.1  & 10.71 \tabularnewline
         97-VA-L & 1280 & -34.99  & 487.58 &  98.81 &  97.3  & 3.57 \tabularnewline
         95-VA-L & 841 & -57.29  & 326.85 &  97.99 &  95.2  & 3.57 \tabularnewline
         90-VA-L & 634 & -67.8  & 281.18 &  96.19 &  90.0  & 3.57 \tabularnewline \bottomrule
        \end{tabular}

}
{$\overline{V}_\text{feas}$ - mean vehicle feasibility in percent, 
$V_{\text{feas}}^{\min}$ - minimum vehicle feasibility in percent,
$Z_\text{feas}$ -  scenario feasibility in percent, $o$ - cost of solution, $\Delta o$ -  cost difference to \gls{abk:FSA} cost, $t$ - computation time in minutes }
\end{table}

The \gls{abk:FSA} solution has very high costs, which amount to more than twice of the median cost over all individual scenarios. 
Changing the $\alpha$-level of the \gls{abk:aSA}, which bounds the minimal scenario feasibility, shows only a small change in costs. Even at a very low $\alpha$-level, requiring feasibility in only 50\% of the scenarios ( $\alpha=0.50$), the costs are only 16.86\% lower than in the \gls{abk:FSA}. However, changing the $\alpha$-level in the \gls{abk:aVA}, which bounds the minimal vehicle feasibility instead, has a very strong influence on costs. Even at a level of $0.99$, allowing 1\% of vehicles being infeasible, the costs are 21.41\% lower for the 99-VA-M compared to the \gls{abk:FSA}. Lowering the level to $0.97$, doubles the reduction to 41.81\%.
The \gls{abk:aVA}-M, shows a stronger cost reduction than the \gls{abk:aVA}-L but is limited at a maximum cost reduction of 50.99\%. The ISA-L forms the lower bound for costs, as by definition it is the lowest cost configuration for which at least one scenario is feasible.

The FSA has the lowest computation time among all approaches.
We see a strong increase in computation time when lowering $\alpha$-levels in the  \gls{abk:aSA}, because allowing some scenarios to remain infeasible,  reduces the possibility for exiting early when evaluating operational feasibility.
In the \gls{abk:aVA} however, we see a reduction of computation time for lower $\alpha$-levels. Since we add scenarios iteratively using the adversarial sampling method, lower feasibility requirements lead to fewer scenarios being included in the optimization.

\begin{figure}[h!]
\FIGURE
{%
    \includegraphics[]{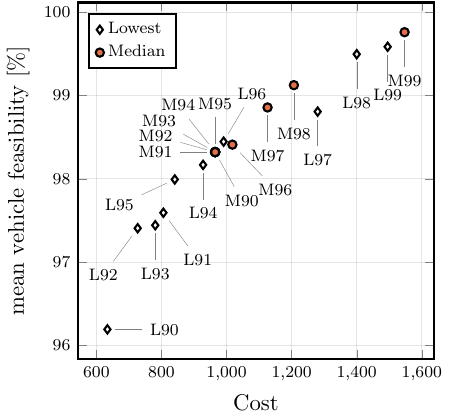}%
    \includegraphics[]{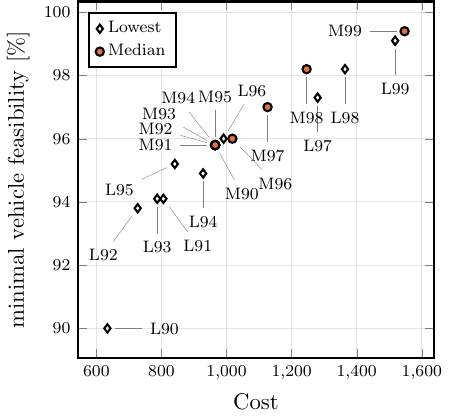}%
}
{Results of \gls{abk:aVA} solution approaches with varying $\alpha$-levels and seed scenarios\label{fig:pareto_benevolent_robust}. 
Labels $M\alpha$ / $L\alpha$ designate seed scenario type medium / lowest and alpha-level.}
{}
\end{figure}

The scenario feasibility is above the desired $\alpha$-levels for the \gls{abk:aSA}. For the \gls{abk:aVA} the scenario feasibility drops to 10.71\% for the 97-VA-L (which ensures 99\% vehicle feasibility), and lowers to 3.57\% for the 97-VA-L. 
The maximal vehicle feasibility is above the chosen $\alpha$-level for all \gls{abk:aVA} approaches. Figure~\ref{fig:pareto_benevolent_robust} shows the relationship between the cost and feasibility for the \gls{abk:aVA}  at varying $\alpha$-levels.  The \gls{abk:aVA}-M cannot reduce the costs below the costs of its seed scenario and is thus unable to reach lower-cost solutions that would be allowed under the given $\alpha$-level. For the upper end of mean feasibilities, the  \gls{abk:aVA}-M variant forms the Pareto frontier, giving slightly better results than the  \gls{abk:aVA}-L.
However, the  \gls{abk:aVA}-L variant dominates at lower vehicle feasibilities,  which the \gls{abk:aVA}-M can not reach.
The \gls{abk:aVA}-L allows a planner much more control, especially when having a larger risk tolerance.

The \gls{abk:aSA}'s ability to accurately select a 
scenario feasibility has little practical relevance in our case study, as the vehicle feasibility remains above 99\% 
even for solutions where 50\% of the scenarios lead to infeasible operational schedules. Accordingly, these scenario infeasibilities affect only a small fraction of the vehicles.
Therefore, we focus on the vehicle infeasibility in the remainder of the study.

\subsubsection{A-Posteriori Validation of Uncertainty Set Approximation}\label{subsec:validation_of_uncertainty_approx}

We obtained solutions using an approximated uncertainty set $\mathcal{Z}$ containing 28 24-hour scenarios of \revnew{vehicles' customer schedules} derived from real-world data. To verify the suitability of this uncertainty approximation, we apply a validation test using six out-of-sample scenario sets denoted as $\mathcal{Z}^1$ to $\mathcal{Z}^6$, each created using the same sampling method that was used to generate $\mathcal{Z}$. We obtain an optimal charging station configuration using the original uncertainty set~$\mathcal{Z}$, and then test this solution for operational feasibility using each of the out-of-sample scenario sets. Specifically, we test if the solution's vehicle feasibility, (for the \gls{abk:aVA}) and scenario feasibility (for the \gls{abk:aSA}) are above the limit set by the $\alpha$ parameter.

\begin{figure}[h!]
    \FIGURE
    {%
        \includegraphics[]{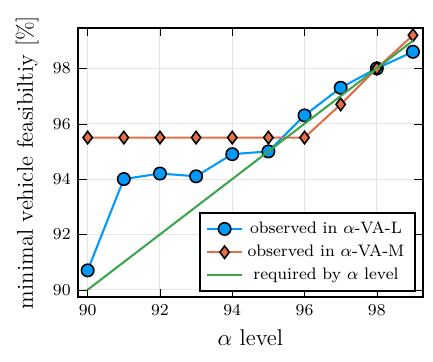}%
        \includegraphics[]{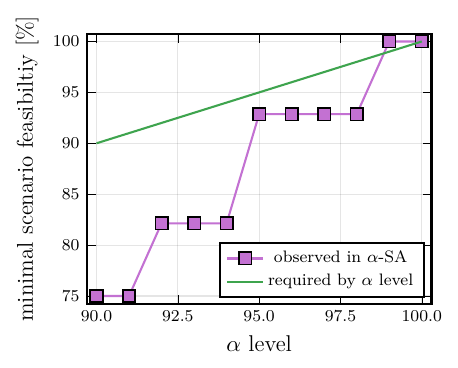}%
    }
    {A-posteriori validation of $\alpha$-VA and  $\alpha$-SA solutions using varying $\alpha$-levels and out-of-sample scenario sets $\mathcal{Z}^1$ to $\mathcal{Z}^6$ \label{fig:aposteriori_trend}.}
    {}
    \end{figure}

Table~\ref{tbl:posteriori_validation} shows the resulting mean and minimum vehicle feasibility when validating the 95-VA-L using the out of sample scenario sets.
While the minimum vehicle feasibilities are slightly lower, none violates the 95\% $\alpha$-level that was selected as feasibility requirement. 
For the 90-SA, however, Table~\ref{tbl:posteriori_validation_sa} shows that the scenario feasibility of all out of sample scenarios violates the 90\% $\alpha$-level.
Figure~\ref{fig:aposteriori_trend} shows the minimal observed vehicle feasibility for the \gls{abk:aVA} and the scenario feasbility for the \gls{abk:aSA} at varying $\alpha$-levels for the same six out-of-sample scenario sets. For the \gls{abk:aVA}, the difference between the allowed and the observed vehicle feasibility is larger
at lower $\alpha$-levels than for higher levels. The \gls{abk:aVA}-M can not generate new solutions below a 96\% $\alpha$-level. Above this value, the minimal a-posteriori vehicle feasibility is lower than the requested $\alpha$-level. The \gls{abk:aVA}-L is exceeding the required $\alpha$-level for all but the 99 $\alpha$-level. For the \gls{abk:aSA} larger $\alpha$ values lead to a smaller difference between observed scenario feasibility and $\alpha$-level. However, only the 99-SA fulfills the required scenario feasibility.

Since the violations of the required vehicle feasibility are minimal, we conclude that scenario set $\mathcal{Z}$ is suitable for creating robust solutions using the \gls{abk:aVA}.
The results further support the preference for the \gls{abk:aVA}-L over the \gls{abk:aVA}-M approach, and the disadvantages of the \gls{abk:aSA} approach.

\begin{table}[h!]
\small

\end{table}

\begin{table}[h!]
    \small

    \begin{minipage}{0.45\textwidth}
        \TABLE
{A-posteriori validation of 95-VA-L solution using out-of-sample scenario sets $\mathcal{Z}^1$ to $\mathcal{Z}^6$.\label{tbl:posteriori_validation}}
{%

    \begin{tabular}{lll}
\toprule
 Scenario Set                   &    $\overline{V_{\text{feas}}}$        & $V_{\text{feas}}^{\min}$     \\\midrule
$\mathcal{Z}$                 &  97.99      &  95.20 \vspace{3pt}  \\
$\mathcal{Z}^1$                 &  98.10      &  95.10   \\
$\mathcal{Z}^2$                 &  98.04      &  95.30   \\
$\mathcal{Z}^3$                 &  98.09      &  96.00   \\
$\mathcal{Z}^4$                 &  98.05      &  95.50   \\
$\mathcal{Z}^5$                 &  98.02      &  95.00   \\
$\mathcal{Z}^6$                 &  98.07      &  95.10   \\
\bottomrule
\end{tabular}
}
{
$\overline{V_{\text{feas}}}$ - mean vehicle feasibility in percent,
$V_{\text{feas}}^{\min}$ - minimum vehicle feasibility in percent 
}
\end{minipage}%
\begin{minipage}{0.45\textwidth}
    \TABLE
    {A-posteriori validation of 90-SA solution using out-of-sample scenario sets $\mathcal{Z}^1$ to $\mathcal{Z}^6$.\label{tbl:posteriori_validation_sa}}
    {%
    
        \begin{tabular}{ll}
    \toprule
     Scenario Set                   &    $\overline{Z_{\text{feas}}}$     \\\midrule
    $\mathcal{Z}$                 &  97.99     \vspace{3pt}  \\
    $\mathcal{Z}^1$                 &  82.14  \\
    $\mathcal{Z}^2$                 &  85.71  \\
    $\mathcal{Z}^3$                 &  78.57  \\
    $\mathcal{Z}^4$                 &  78.57  \\
    $\mathcal{Z}^5$                 &  75.00  \\
    $\mathcal{Z}^6$                 &  89.29  \\
    \bottomrule
    \end{tabular}
    }
    {
    $Z_{\text{feas}}$ - scenario feasibility in percent 
    }
\end{minipage}%
    \end{table}

\subsubsection{Sensitivity to Technological Developments}\label{section:computational_study:solution_analysis}\label{sec:comp:tech_change}

We analyze the sensitivity of cost and mean vehicle feasibility to battery and charging parameter changes using the  \gls{abk:aVA}-L algorithm with an $\alpha$-level of 95. To this end, we vary the battery size and the charging speed between -50\% and +50\% relative to the baseline configuration.
Changes in charging technology might also allow for depot charging between taxi shifts for the full fleet. Accordingly, we model depot charging of amount $q^d$ by reducing the SOC requirement at the end of the shift $q^\text{e}_v$ by $q^d$, i.e., $0.5 = q^\text{b}_v = q^\text{e}_v + q^d$ holds. 
Since no depot charging $(q^d = 0)$ is considered in the baseline scenario, we illustrate the relative variation of the depot charging parameter as absolute $q^d$ change, e.g., a $+20\%$ variation leads to $(q^d = 0.2)$, which in turn represents 20\% of the total battery capacity being charged at the depot.

\begin{figure}[h!]
\FIGURE
{
    \includegraphics[]{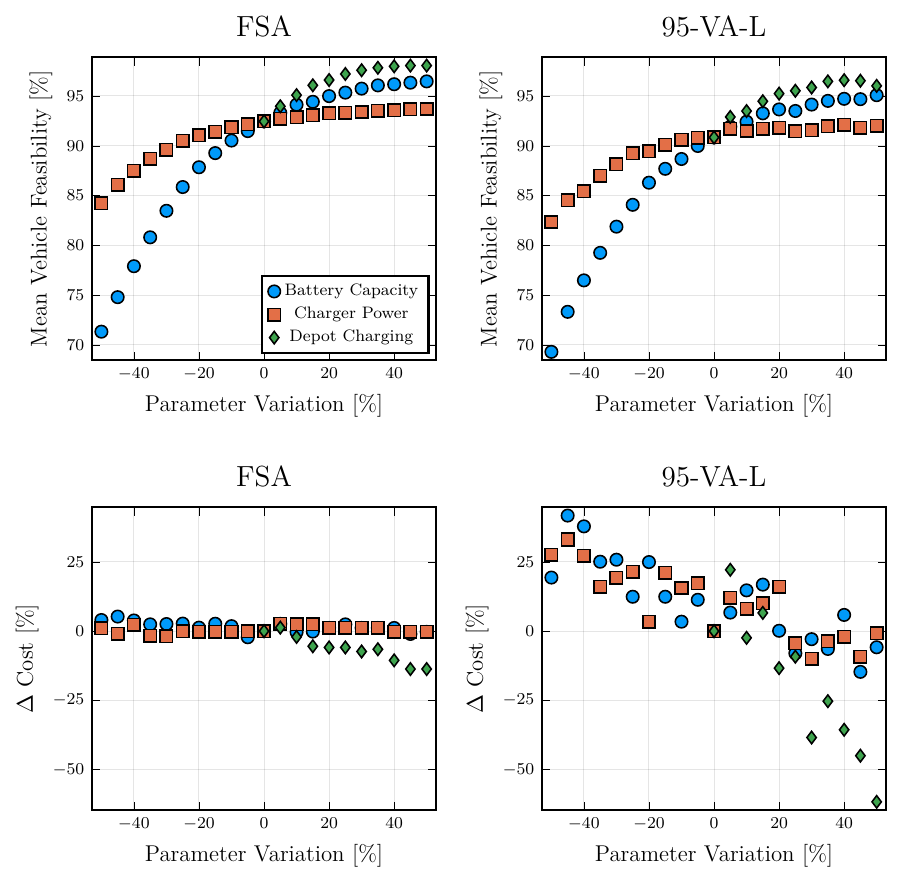}
    }
{Sensitivity of mean vehicle feasibility and cost of the solution for both FSA and 95-VA-L to variation of technological parameters. \label{fig:rc_sensitivity}}
{}
\end{figure}

We show the mean vehicle feasibility on the master instance and relative cost change for modifying these parameters ceteris paribus in Figure~\ref{fig:rc_sensitivity}. 
The vehicle infeasibility shown for the FSA originates entirely from the vehicles removed for the parameter-specific instance since the FSA otherwise ensures full vehicle feasibility on all scenarios considered.
The feasibility of the 95-VA-L largely follows the same trend, with each feasibility being slightly lower than in the FSA (1.5 pp for the baseline case), showing that the parameter-specific feasibility is also dominant in the 95-VA-L. 

The feasibility sensitivity shows that the battery capacity changes have a larger impact than the charging speed, with both being dominated by depot charging.   However, varying only the state of charge required at the end of the shift as a proxy for depot charging overestimates the impact of depot charging since a detour requirement to a potential depot and additional investment requirements are omitted. Nevertheless, optimizing the internal processes by allowing depot charging during the vehicle handover will likely be very effective. The relative change of costs shows that deteriorating the technological parameters leads to increased costs. We note that the costs are not directly comparable since the number of feasible vehicles also changes. We have higher costs with a lower battery capacity and can also serve fewer vehicles than in the baseline configuration.
A reduction of battery capacity by 50\% results in a cost increase of more than 30\% for the 95-VA-L, while the FSA shows little change in costs. However, the FSA's baseline solution costs are with 1846 already more than twice the cost than the 822 of the 95-VA-L baseline. Due to the cost difference, the FSA solution has the slack available to cover the higher demand caused by smaller batteries or slower charging speeds without additional investments.

\revnew{
\subsubsection{Variable Charging Station Sizing}\label{section:computational_study:variable-sizing}

We will now consider the variable-size approach introduced in Section~\ref{sec:variable_sizing}. We use the instance with 60 potential charging station locations and first analyze the results of the deterministic approach for increasing vehicle numbers. We then show the results of the robust 95-VA-L approach on the $|N_s| = 60, |V| = 1,000$ instance.

Figure~\ref{fig:variable_runtime_duration} shows a box-whisker plot of computation times of the variable-size approach for increasing numbers of vehicles.  The factor by which the average computation time of the variable-size approach increases compared to the static-size approach is shown on the right axis. The mean computation time of the variable-size approaches with 1,000 vehicles is at 2,608 seconds, which is 4.6 times the runtime of the static-size approach. Further, the relative difference to the static-size case increases with increasing vehicle numbers.
The robust setting shows comparable results when using a variable-size 95-VA-L approach. Its computation time increases by a factor of 4.5 compared to the static-size variant.

\begin{figure}[h!]
    \FIGURE
    {
        \includegraphics[]{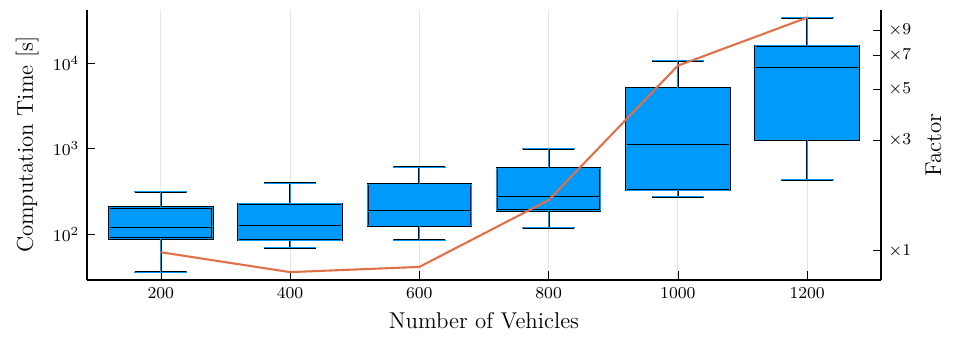}
        }
    { The blue box plots give the computation time (log scale) of the variable-size algorithm for instances with $|\mathcal{N}^s| = 60$ and variable number of vehicles (left axis). The orange line gives the factor by which the average computation time of the variable-size approach increases compared to the static-size approach (right axis).
    \label{fig:variable_runtime_duration}}
    {}
\end{figure}

Interestingly, the vast majority of instances have optimal solutions with only two charging points for each selected location. For the instances where the optimal solution contains locations with two and four charging points, a typical pattern is that the sites with four charging points are located at the border between the inner- and the outer-city area, profiting from the lower location costs of the outer-city area but still close to the majority of the demand stemming from the inner-city area.

Considering the cost advantages of larger charging stations, the long lifetime of public infrastructure and
increasing future demand, the selection of four charging points per charging station as preferred by the local utility company remains reasonable.

}

\revnew{

\subsection{Discussion of the applicability of our results for the real-world setting}

Some of the assumptions we have made to model the problem may seem limiting from a practical perspective. In the following, we present dedicated experiments and studies to investigate these limitations.

\paragraph{Assessment of perfect information solution in an online setting}\label{sec:simulation}
Although our approach accounts for uncertainty in vehicles' customer schedules by leveraging robust optimization, the resulting planning problem still simplifies the operational problem: in practice, drivers must make online charging decisions without knowing about future customer demand or charging station availability. However, accounting for such a dynamic in the strategic planning problem at hand is hardly tractable and would require to combine a simulation with an optimization. To analyze the error that our model inherits with this simplification, we report on a simulation study that we have undertaken to evaluate the quality of the charging station decision obtained with our offline approach for a dynamic online setting. The simulation focuses on the impact of online charging decisions and does not consider vehicles competing for customers, i.e., the customers of vehicle $v$ are replayed from its historical customer schedule $C_v$. If the battery of a vehicle reaches an SOC of zero while being in operation, we count it as infeasible. 
The decision space in the simulation follows the options available in the time-space network of the original operational problem. We assume that drivers know their current SOC and can access real-time information about the locations and availability status, i.e., the number of free charging points, at charging stations. Accordingly, we assume that drivers use this information to make a charging decision after each customer drop-off based on two SOC thresholds: (i) If the SOC is below the minimum charge threshold  $T^-$, the driver must make a detour to charge. To that end, the driver evaluates the availability of the three nearest stations and drives to the closest available. If none are available, the driver drives to the closest station and expects to wait. (ii) If the SOC is below the opportunity threshold  $T^o$, the driver will drive to a charging station if it is available and reachable within 5 minutes. 

Upon arrival at a charging station, drivers will charge until the SOC reaches the maximum charge threshold $T^+$. If all charging points are occupied at the time of arrival, drivers will wait until a charging point becomes available. Drivers must terminate waiting or charging if a customer request arrives. We sample thresholds $T^-$ and $T^o$   for each driver to account for drivers' varying range anxieties. In accordance with \cite{trippe_mobility_2015} we set $T^-$ to follow a normal distribution $T^- \sim N(0.30, 0.05)$.  When drivers have completed 75\% of their shift, we raise the threshold to $T^- \sim N(0.60, 0.05)$ because drivers anticipate that they must return the vehicle at 50\% SOC. We assume that drivers stop charging once reaching 90\% SOC, $T^+ \sim N(0.90, 0.05)$, and that they opportunistically select charging stations below 60\% SOC, $T^o \sim N(0.60, 0.05)$. 

We use the scenarios and their optimal charging station locations of the instance $|N_s| = 60, |V| = 1,000$  for our evaluation. For each scenario, we run 100 trials of our simulation and derive averages. Across all scenarios, the median number of infeasible vehicles is 15.8\%. This shows that the results of our optimization method perform reasonably well when evaluated in an online setting.

Figure~\ref{fig:plot_simulation} shows the distribution of the percentage of infeasible vehicles for the static-size case (Section~\ref{sec:solution_method:stategic_problem}) and the variable-size case (Section~\ref{sec:variable_sizing}). We can see that the median infeasibility rate of the variable-size case is 33\%. The infeasibility is more than twice as high as in the static-size case. The simulation in an online setting shows that the lower cost of the variable-size approach compared to the fixed-size approach comes at a substantial risk of operational infeasibility. The reliability requirements of commercial fleets, the lower marginal costs of larger charging stations, and the growth of charging demand in the future are speaking in favor of four charging points.

\begin{figure}[h!]
    \FIGURE
    {
        \includegraphics[]{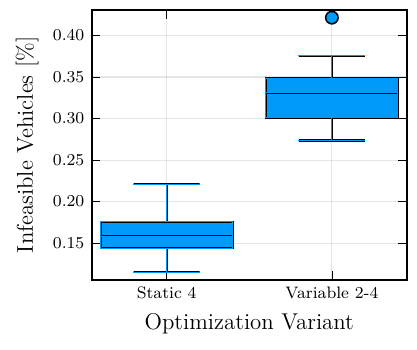}
        }
    {Percentage of infeasible vehicles when the online solution is applied in an online setting for $|\mathcal{N}^s| = 60$ and $|V| = 1,000$. \label{fig:plot_simulation}}
    {}
\end{figure}

\paragraph{Plausibility of planned charging detours}

In our approach, we consider only the cost of constructing charging stations in the objective function.
The distances of detours are only indirectly penalized via the additional energy consumption, potentially leading to excessive detours. In the following, we use the optimal solution, i.e., opened stations and scheduled charging operations, of the first scenario of the instance with $|V|=1,000, |N^s|=60$ to show that the optimal solutions do not contain large detours.

We first determine the distance increase incurred through charging detours by comparing the distance driven when including the detour to the direct distance between the customers. Since we have computed the fastest routes on real road network data for our distance matrix, for a few instances, the detour leads to shorter distances than the direct route, and thus negative distance increases. For those cases, we set the distance increase to zero.
 On average, the distance driven in a shift increases by only 2.9\% when including charging detours. Figure~\ref{fig:inc_dist_detour_shift} shows a histogram of the distance increase for all taxi shifts in the scenario. 

\begin{figure}[h!] \centering
    \includegraphics[]{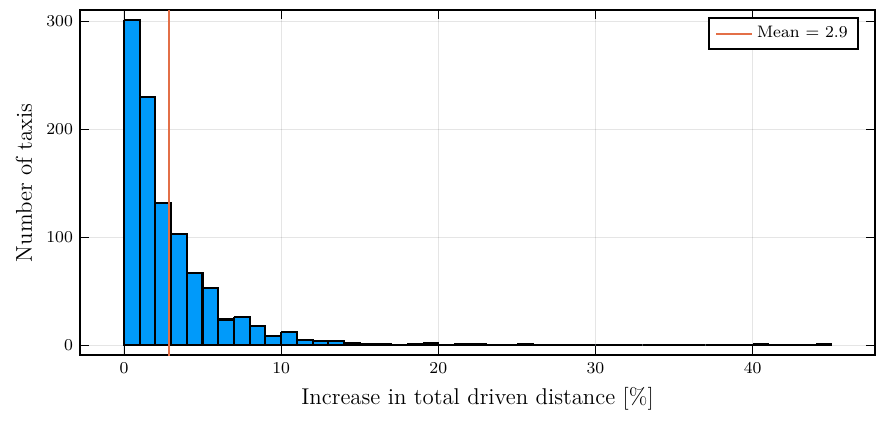}
    \caption{The distance increase resulting from charging relative to the total driven distance of each taxi} \label{fig:inc_dist_detour_shift}
\end{figure}

Next, we validate for each scheduled charging detour if the assigned charging station is the nearest reachable. We find that in 55\% of the cases, the assigned charging station is the nearest.
If we assigned a more distant charging station, it is, on average, only 1.3 km further away than the nearest one.
}

\section{Conclusion}\label{sec:conclusion}

This paper proposes a novel robust solution method to design
optimal charging infrastructure for a sizeable ride-hailing fleet. Our solution method takes a novel perspective by working with real-world customer schedules, uncertain demand, capacitated charging stations, non-linear charging functions and making individual charging and detour decisions for every vehicle in the operational problem. 
We develop an exact algorithm based on cutting planes and branch-and-price, embedded into a robust optimization framework.
We introduce a \gls{abk:FSA}, an \gls{abk:aSA}, and an \gls{abk:aVA} algorithm. 
The latter two algorithms allow a planner to control the risk tolerance by adjusting the allowed infeasibility through an $\alpha$ parameter: In the \gls{abk:aSA}, we relax the requirement that all scenarios have to be feasible, while in the \gls{abk:aVA}, we allow for a definable fraction of infeasible vehicles within the scenarios using an adversarial sampling procedure.

We perform a computational study using real-world data of a Munich. Our computational study shows the efficacy of our algorithms, supporting scenarios containing more than 1,000 individual vehicles serving approximately 11,000 customers.  
In the analysis of the robust variants, we show that the $\alpha$-VA-L variant is most suitable for the problem at hand. The $\alpha$ parameter can be set depending on the risk tolerance of the planner, with the $\alpha  = .95$ being a good trade-off between feasbility and cost. Allowing even a 1\% feasibility tolerance on individual vehicles leads to significantly lower costs. In an a-posteriori feasibility analysis, we confirm that feasibility requirements are met by validating the solutions through out-of-sample scenario sets.
Our sensitivity analysis of technological parameters shows that allowing for depot charging before or during the vehicle hand-over has a more significant
effect on cost than focusing on battery capacity and charging speed. The battery capacity has a slightly stronger influence on cost and
feasibility than the available charging speed. 
\revnew{
    In this study, we operate on fixed customer tours as seen in historical data. While this allows us to show that electrification is possible for existing driving patterns, rerouting 
    customers constitutes a promising direction for future research.  Using a modified version of our branch-and-price approach, we have developed a prototype that allows customers to be reassigned to vehicles. However, we can only solve small instances, which necessitates the investigation of alternative solution approaches. 
}



\ACKNOWLEDGMENT{%
We thank Michael Wittmann from the Institute of Automotive Technology at the Technical University of Munich
for providing and assisting with the real world taxi data set. 
The authors gratefully acknowledge the computational and data resources provided by the Leibniz Supercomputing Centre (www.lrz.de).
This work has been partially funded by the Deutsche Forschungsgemeinschaft (DFG, German Research Foundation) - Project Number 277991500

}

%
\begin{APPENDIX}{Notation}

\begin{longtable}{llp{0.6\textwidth}}
    \midrule
    {\small Problem Description:} &  & \\   
    & $\mathcal{V}$ & Set of vehicles\\
    & $\mathcal{C}_v$ & Ordered Set of customers per vehicle \\
   
    &${\cal N}^s$ & Set of candidate locations for charging stations \\
    &${\cal N}^c$ & Set of customer pickup and drop-off locations \\
    &$\mathcal{S}$ &  Set of opened charging stations\\
    & $\mathcal{T}$ & Set of discrete times \\
    & $\mathcal{Z}$ & Uncertainty set of \revnew{vehicle routes} \\

    & $\cal {G} = ({\cal N},{\cal A})$  & Problem graph with node and arc sets \\

    &$d_{ij}$ & Distance of arc \\
    &$\tau_{ij}$ & Travel time of arc \\
    & $q_{v,t}$ & SOC of vehicle $v$ in time $t$ \\
    &$e_{ij}$ & Energy consumption of arc \\
    &$L^o_c$ & Customer pickup location\\
    &$L^d_c$ & Customer drop-off location\\
    &$T^o_c$ & Customer pickup start time  \\
    &$T^d_c$ & Customer drop-off arrival time  \\
    & $\gamma(q,l)$ & Charging function of (duration $l$ and SOC at time of charging $q$) \\
    & $q^b_v$ & SOC of vehicle $v$ at beginning of its shift \\
    & $q^\text{e}_v$ & SOC requirement of vehicle $v$ at end of its shift \\
    & $[q^{\min}, q^{\max}]$ & Allowed battery range \\

    & $o_s$ & Cost of opening charging station\\
    & $n_s$ & Number of chargers at charging site \\  \midrule

    {\small Cutting Plane Problem:} &  & \\

    & $\mathcal{C}$ & Cuts in CCP \\
    & $\mathcal{I}_c$ & Cover set of charging stations inducing infeasibility in cut $c$ \\
    & $\mathcal{C}$ & Set of active cover sets \\
    & $\mathcal{O}$ & Set of improved cover sets \\
    & $x_s$ & Opening decision for charging station s\\  \midrule

    {\small Operational Problem:} &  & \\

    & $\mathcal{H}_v$ &  Set of customer vertices \\
    & $\mathcal{R}_v$ &  Set of shift arcs representing vehicle movements \\
    & $\mathcal{Q}_v$ &  Set of charging vertices \\
    & $\mathcal{B}_v$ &  Set of charging and movement arcs \\
    & $\mathcal{W}_v$ &  Set of holdover arcs at charging station \\
    & $\mathcal{P}_v$ &  Set of paths for vehicle $v$ \\   
    
    & ${\cal {G}}_v = ({\cal N}_v,{\cal A}_v) $ & Time expanded graph for vehicle $v$ \\

    &$T^t_{cs}$ & Travel time between customer drop-off and charging station \\
    &$T^t_{sc}$ & Travel time between charging station and customer pick-up \\

    &$k_{pst}$ & Parameter showing if path $p$ uses charging station $s$ in time $t$ \\
    &$\lambda_v$ & Dual of vehicle convexity constraint \\
    &$\rho_{st}$ & Dual of capacity constraint \\
    &$c_p$ & Reduced cost of $x_p$ \\

    & $\pi_{ij}$ &  Charging costs of arc ij\\
    & $z_{ij}$ &  Charging duration before arc ij\\
    & $r_{ij}$ & Range costs of arc ij\\
    & $L_i = ( C_i, R_i)$ & Label $i$ with accumulated charging costs $C_i$ and SOC $R_i$ \\

    & $y_{vp}$ & Path $p$ selected for vehicle $v$ \\
    & $i_v$ & Dummy variable for vehicle $v$ \\

    {\small Robust Framework:} &  & \\

    & $\Omega$  & Full feasibility set \\
    & $\Phi$  & Relaxed feasibility set \\
    & $z^*$ & Scenario candidate\\
    & $z'$ & Partial scenario derived from $z^*$ \\
    & $V_z$ & Vehicles contained in scenario $z$ \\
    & $V_z^f$ & Vehicle feasibility of scenario $z$ \\
    & $V_\text{feas}^{\min}$ & Minimal vehicle feasibility \\
    & $\overline{V}_\text{feas}$ & Mean vehicle feasibility \\
    & $ Z_\text{feas}$ & Scenario feasibility \\
    \bottomrule
\end{longtable}

\end{APPENDIX}
%
%


\bibliographystyle{informs2014trsc} 
\bibliography{etaxi2019.bib,main.bib,oscm.bib} 


\end{document}